\documentclass[11pt]{amsart}

\usepackage{subcaption}

\usepackage[margin=1.25in,
           top=1.25in,
           bottom=1.25in]{geometry}
\usepackage{amssymb}
\usepackage{amscd}
\usepackage{mathtools}
\usepackage[all]{xy}
\usepackage{bbm}
\usepackage{mathrsfs}
\usepackage{enumerate}
\usepackage{tikz}
\usepackage{stmaryrd}
\usepackage{etoolbox}
\usepackage{array}
\usepackage{hyperref}
\usepackage[lite]{amsrefs}
\usepackage{comment}

\newcommand{\N}{\mathbb{N}}
\newcommand{\Z}{\mathbb{Z}}

\newcommand{\R}{\mathbb{R}}

\newcommand{\eps}{\varepsilon}

\newcommand{\del}{\nabla}
\newcommand{\ind}{\operatorname{ind}}

\newcommand{\lap}{\Delta}

\newcommand{\bd}{\partial}
\newcommand{\cl}{\overline}
\newcommand{\ins}{\operatorname{int}}

\newcommand{\tr}{\operatorname{tr}}

\newcommand{\la}{\langle}
\newcommand{\ra}{\rangle}

\newcommand{\dist}{\operatorname{dist}}

\renewcommand{\div}{\operatorname{div}}

\newcommand{\grad}{\del}
\newcommand{\vol}{\operatorname{Vol}}

\newcommand{\area}{\operatorname{Area}}

\newcommand{\ric}{\operatorname{Ric}}

\newcommand{\mb}{\mathbb}
\newcommand{\mc}{\mathcal}

\newcommand{\mk}{\mathfrak}

\newcommand{\wha}{\widehat}

\newcommand{\Area}{\mathrm{Area}}
\newcommand{\Id}{\mathrm{Id}}
\newcommand{\diag}{\mathrm{diag}}

\newcommand{\Span}{\mathrm{span}}

\newcommand{\im}{\mathrm{Im}}
\newcommand{\Rm}{\mathrm{Rm}}

\newcommand{\codim}{\operatorname{codim}}

\newcommand{\rom}[1]{\expandafter\romannumeral #1}
\newcommand{\Rom}[1]{\uppercase\expandafter{\romannumeral #1}}

\DeclareMathOperator{\Ric}{Ric}

\DeclareMathOperator{\Diff}{Diff}

\DeclareMathOperator{\interior}{int}

\DeclareMathOperator{\Ein}{Ein}

\newcommand{\blue}{\color{blue}}

\theoremstyle{plain}
\newtheorem{theorem}{Theorem}[section]

\newtheorem{prop}[theorem]{Proposition}
\newtheorem{lem}[theorem]{Lemma}
\newtheorem{conj}[theorem]{Conjecture}
\newtheorem{claim}[theorem]{Claim}

\theoremstyle{definition}

\newtheorem{rem}[theorem]{Remark}

\newtheorem{proposition}[theorem]{Proposition}
\newtheorem{lemma}[theorem]{Lemma}

\theoremstyle{definition}
\newtheorem{definition}[theorem]{Definition}
\theoremstyle{remark}
\newtheorem{remark}[theorem]{Remark}

\numberwithin{equation}{section}

\begin{document}

\title[PE Liouville and Yamabe Invariant]{A Poincare-Einstein Liouville Theorem and the Yamabe Invariant of Lens Spaces}
\author{Liam Mazurowski}
\address{Department of Mathematics, Southern Connecticut State University, New Haven, CT 06515}
\email{mazurowskil1@southernct.edu}

\author{Tongrui Wang}
\address{School of Mathematical Sciences, Shanghai Jiao Tong University, Minhang District, Shanghai 200240, China}
\email{wangtongrui@sjtu.edu.cn}

\author{Xuan Yao}
\address{Department of Mathematics, University of Chicago, Chicago, IL 60637}
\email{dorisxyao@uchicago.edu}

\author{Xin Zhou}
\address{Department of Mathematics, Cornell University, Ithaca, New York, 14853, United States}
\email{xinzhou@cornell.edu}

\begin{abstract}
   Schoen conjectured that the Yamabe invariant of every lens space is achieved by the round metric. A proof of this conjecture is presented for all lens spaces, together with rigidity: equality is only achieved for the round class. The key step is to show that every $\mathbb Z_p$-invariant Yamabe-positive conformal class on $S^3$ admits a $\mathbb Z_p$-equivariant Poincaré-Einstein filling. The distance function from a fixed point of the filling then gives an invariant test function with Yamabe quotient at most $\sigma(S^3)$, which descends to the lens space.
 The proof of the main theorem was generated by AI (ChatGPT Astra, September 2026); this is mainly a learning note of the AI generated proof.

\end{abstract}

\maketitle

\section{Introduction}  The Yamabe invariant  of a smooth manifold is an important invariant with close links to topology (cf. \cite{lee1987yamabe}).  Nevertheless, computing the Yamabe invariant of specific manifolds has proven to be a difficult problem. Focusing on the three dimensional case, it is known that $\sigma(M)$ is achieved by the hyperbolic metric on hyperbolic manifolds \cite{anderson2006canonical} and by the flat metric on tori.  In positive curvature, the value of $\sigma(S^3)$ is known, and Bray-Neves \cite{bray2004classification} computed the value of $\sigma(\mathbb{RP}^3)$ by exploiting a connection with the Schwarzschild manifold.  Schoen \cite{schoen1987variational} conjectured that the Yamabe invariant of any spherical space form should be achieved by the round metric. 

\begin{conj}[Schoen \cite{schoen1987variational}] \label{schoen-conjecture} For any spherical space form $M = S^3/\Gamma$ one has $\sigma(M) = p^{-2/3} \sigma(S^3)$ where $p = \#\Gamma$.  
\end{conj}

\subsection*{AI Disclosure} 
The purpose of these notes is to explain an AI generated proof of Conjecture \ref{schoen-conjecture} for Lens spaces plus a rigidity statement (Theorem \ref{sigma-lens}), obtained using Chat GPT Astra in September 2026. We will focus on explaining the novel ideas in the proof. No AI generated text was included. Similar results were presented independently in the recent arXiv preprint \cite{loren}. 

\subsection{Poincar\'e-Einstein Metrics and the Yamabe Invariant}
It has been known for some time that there is a connection between the Yamabe invariant and Poincar\'e-Einstein fillings.  To the authors' knowledge, such a connection was first found by Dutta-Javaheri \cite{dutta2010rigidity} and then further developed by Li-Qing-Shi \cite{li2017gap}.  Suppose $h$ is a Yamabe positive metric on $S^3$ and $g_+$ is a Poincar\'e-Einstein filling on $B^4$.  Fix a point $p_0\in B^4$ and define $r = d_{g_+}(\cdot,p_0)$. Also, let $x$ be a geodesic defining function in a collar neighborhood of the boundary. Define the function $u = r + \log x$ on $B^4$. In \cite{li2017gap}, they show that $u$ extends to a Lipschitz function $u_0$ on $S^3$. Then, by a delicate analysis of the geometry of the level sets of $r$, they show that the function $\psi_0 := e^{-u_0/2}$ satisfies the following pair of inequalities
\begin{gather}
    \label{in1} \int_{S^3} \psi_0^6\, dv_h \le \frac{\omega_3}{8},\\
    \label{in2} \int_{S^3} \vert \grad^h \psi_0\vert^2 + \frac 1 8 R_h \psi_0^2\, dv_h \le 3 \int_{S^3}\psi_0^6 \, dv_h,
\end{gather}
where $\omega_3=2\pi^2$ is the volume of the unit round $S^3$.
Combining these, it then follows that $Y(S^3,[h])\le \sigma(
S^3)$. 

Now for $p\in \{2,3,\dots\}$ and coprime $1\leq q<p$, consider $\mk r_{p,q}\colon (z,w)\mapsto (e^{2\pi i/p} z,e^{2\pi iq/p}w)$ so that $S^3/\la \mk r_{p,q}\ra$ is the Lens space $L(p,q)$.  
Assume that $h$ is a metric on $S^3$ invariant under the action of $\mk r_{p,q}$.  The first key idea is that, given the existence of an {\it equivariant} Poincar\'e-Einstein filling, the above argument can all be made equivariant. Indeed, suppose $\phi$ is a diffeomorphism of $\cl B^4$ with $\phi \equiv \mk r_{p,q}$ on $S^3$ and that $g_+$ is a Poincar\'e-Einstein filling of $h$ with $\phi^* g_+ = g_+$. Such a diffeomorphism must have a fixed point in $B^4$, and so we can choose $p_0$ to satisfy $\phi(p_0) = p_0$. Then the distance function $r$ and the defining function $x$ will both be invariant under $\phi$, and it follows that the function $\psi_0$ on the boundary is invariant under $\mk r_{p,q}$. The inequalities \eqref{in1} and \eqref{in2} then descend to the quotient $L(p,q)$ and it follows that $Y(L(p,q),[h]) \le p^{-2/3} \sigma(S^3)$. Thus Schoen's conjecture is reduced to the existence of an equivariant Poincar\'e-Einstein filling. 

\subsection{Existence of Poincar\'e-Einstein Fillings} There is an extensive literature on the Poincare-Einstein fill in problem.  We refer the reader to the work of Chang-Ge \cite{chang2025problem} for more background on this topic. Temporarily putting equivariance to the side, let $h$ be a Yamabe positive metric on $S^3$. By work of Anderson \cite{anderson2010structure} and Chang-Ge \cite{chang2025problem}, the existence of a Poincare-Einstein filling $g_+$ for $h$ on the ball $B^4$ can be reduced to proving the following Liouville-type theorem. 

\begin{theorem}
Let $X^4$ be a smooth, connected $4$-manifold with boundary.  Assume that one boundary component of $X$ is diffeomorphic to $\R^3$ and let $\delta$ be a Euclidean metric on this boundary component.  Assume that $g_+$ is a complete Einstein metric on the interior of $X$ satisfying $\ric(g_+) = -3g_+$.  Assume that for each closed ball $K\subset \R^3$ there is a neighborhood $U$ of $K$ in $X$ and a smooth defining function $\rho$ on $U$ such that $\rho^2 g_+$ is a smooth metric on $U$ which restricts to $\delta$ on $U\cap \R^3$.  Then $g_+$ is isometric to the hyperbolic metric.   
\end{theorem}

\begin{remark}
{
A similar Liouville theorem was established by Chang–Ge \cite[Theorem II]{chang2025problem} under additional hypotheses. Chang–Ge have informed us that they have already strengthened their Liouville theorem and compactness theorem by removing these additional hypotheses in a revised manuscript. We thank them for sharing this development with us.}
\end{remark}

The main strategy of their proof is similar to the above. They obtain a function $\psi_0$ on the Euclidean boundary component and show that it satisfies inequalities analogous to \eqref{in1} and \eqref{in2}. The Euclidean Sobolev inequality then forces these to be equalities. In turn, this forces equality in the Bishop-Gromov volume comparison and this implies that $g_+$ must be hyperbolic. Their proof again relies on a delicate analysis of the level sets of $r$, which is complicated by the possible presence of the cut locus. This accounts for the extra hypotheses they must impose. 

The second key idea is that inequalities \eqref{in1} and \eqref{in2} can both be obtained by analyzing open 4-dimensional regions rather than 3-dimensional hypersurfaces.  This substantially simplifies the analysis that must be carried out. First, for each compact $K\subset \R^3$, there is a  comparison 
\[
\frac{1}{3}\int_{K} \psi_0^6 \, dv_\delta \le \lim_{R\to \infty} e^{-3R} \vol_{g_+} B(p_0,R)
\]
which yields \eqref{in1} when combined with the Bishop-Gromov volume comparison theorem. Second, the distributional Laplace comparison theorem 
\[
\lap_{g_+} r \le 3\coth r
\]
can be passed directly to a distributional inequality for the function $\psi_0$ on the boundary which yields \eqref{in2}. 

\subsection{Organization} The remainder of these notes is organized as follows. In Part 1, we explain the proof of the Poincar\'e-Einstein Liouville theorem. Then in Part 2, we explain how to use the Poincar\'e-Einstein Liouville theorem in conjunction with results of Anderson \cite{anderson2010structure} and Chang-Ge \cite{chang2025problem} to obtain the existence of equivariant Poincar\'e-Einstein fillings. Finally, in Part 3 we explain how this implies Schoen's conjecture for Lens spaces. As part of this step, we give a new simplified proof of the Li-Qing-Shi \cite{li2017gap} result that reflects the new ideas in the proof of the Poincar\'e-Einstein Liouville theorem. This argument is our own and was not suggested by the AI. 

\part{A Poincar{\'e}-Einstein Liouville Theorem}
\section{Statement and outline} 

\begin{theorem}\label{Thm: Liouville theorem for round filling}
Let $X^4$ be a smooth, connected $4$-manifold with boundary.  Assume that one boundary component of $X$ is diffeomorphic to $\R^3$ and let $\delta$ be a Euclidean metric on this boundary component.  Assume that $g_+$ is a complete Einstein metric on the interior of $X$ satisfying $\ric(g_+) = -3g_+$.  Assume that for each closed ball $K\subset \R^3$ there is a neighborhood $U$ of $K$ in $X$ and a smooth defining function $\rho$ on $U$ such that $\rho^2 g_+$ is a smooth metric on $U$ which restricts to $\delta$ on $U\cap \R^3$.  Then $g_+$ is isometric to the hyperbolic metric.   
\end{theorem}

Here a smooth function $\rho$ on $U$ is called a smooth defining function on $U$ provided the following conditions hold: 
\begin{itemize}
\item[(i)] $\rho > 0$ in $\ins(X) \cap U$ and $\rho = 0$ on $\bd X \cap U$,
\item[(ii)] $d\rho \neq 0$ on $\bd X \cap U$,
\item[(iii)] $\rho^2 g_+$ extends to a smooth metric on all of $U$ including $\bd X \cap U$. 
\end{itemize}
By the Einstein condition and the conformal change formula of Ricci curvature, every defining function $\rho$ on $U$ must actually satisfy $\vert \grad^{\rho^2 g_+} \rho\vert = 1$ on $\bd X \cap U$. 

Under the conditions of the above theorem, given any compact set $K\subset \R^3$, we can always find a neighborhood $U$ of $K$ in $X$ in which there is a smooth {\it geodesic} defining function $x$ which induces the Euclidean metric on $U \cap \R^3$.  This means that $x$ is a smooth function such that $g = x^2 g_+$ is a smooth metric on $U$ and the restriction of $g$ to $U\cap \R^3$ is $\delta$ and $\vert \grad^g x \vert = 1$ in $U$. To see this, write $g = x^2 g_+$ as a metric conformal to $\rho^2 g_+$. Then the condition for $x$ to be a geodesic defining function inducing $\delta$ is a 1st order PDE that can be solved in a neighborhood of $K$ using the method of characteristics. This also implies the uniqueness of $x$ on any overlaps since $x$ is obtained by solving an ODE along characteristics emerging from the boundary.  

\subsection{Outline} 
The core strategy of the proof is the same as the proof of Theorem II in Chang-Ge \cite{chang2025problem}. 
Fix a point $p\in \ins(X)$ and define $r(q) = d_{g_+}(q,p)$.   We consider the function 
\[
u(x,z) = r(x,z) + \log(x) 
\]
and show that it induces a well-defined function $u_0$ on $\R^3$.  We then set $\psi = e^{-u/2}$ and $\psi_0 = e^{-u_0/2}$ and aim to show that $\psi_0$ is a minimizer for the Euclidean $L^2$ Sobolev inequality.  Once this is established, it is shown to force equality in the Bishop-Gromov volume comparison theorem for $g_+$ and this implies that $g_+$ is hyperbolic. 

There are two primary simplifications compared with the Chang-Ge proof.  First, we obtain the $L^6$ estimate for $\psi_0$ at infinity by comparing with the volume of regions rather than the area of spheres.  This eliminates the need to make a detailed study of the geometry of geodesic spheres.  Second, we are able to directly pass the distributional Laplace comparison theorem for $r$ on a local collared neighborhood of the boundary to a distributional inequality for $\psi_0$ on the boundary.  Again this eliminates the need to make a delicate analysis of the geometry of distance spheres across the cut locus, and bypasses the need for hypotheses (2) and (3) of \cite[Theorem II]{chang2025problem}.

\section{Defining the Boundary Function} 
\label{section:defn}

Fix a point $p\in \ins(X)$ and define $r(q) = d_{g_+}(q,p)$. Let $K$ be a closed ball in $\R^3$. Let $U$ be a neighborhood of $K$ on which there is a local geodesic defining function and assume that $\eps > 0$ is small enough that 
\[
N = \{(x,z): x\in (0,\eps),\ z\in K\}
\]
is contained in $U$. 
Define the function 
\[
t = -\log x
\]
on $N$.  Then the metric $g_+$ can be written in the form 
\[
g_+ = dt^2 + e^{2t} h_{e^{-t}} = x^{-2}(dx^2 + h_x)
\]
where $h_x = \delta + O(x^3)$. Corresponding expansions also hold for the derivatives of $h_x$. 

Define the function 
\[
u(x,z) = r(x,z) + \log x. 
\]
The goal of this section is to show that $u$ satisfies 
\[
\vert \bd_z u \vert \le C, \quad 0\le \bd_x u \le C x 
\]
almost-everywhere in $N$. These constants are allowed to depend on $K$. This will then imply that $u$ extends to a Lipschitz function $u_0$ on $K$.  Since the geodesic defining function is well-defined on overlaps, we can then take an exhaustion of $\R^3$ by compact sets and obtain a locally Lipschitz function $u_0$ on all of $\R^3$. 

\subsection{Estimates on Geodesics}
To prove these bounds we need to understand the derivatives of $r$, and this will involve studying minimizing geodesics connecting $p$ to a point $q\in N$.  Recall that 
\[
g_+ = dt^2 + e^{2t} h_{e^{-t}}.
\]
Consider a unit speed $g_+$ geodesic $\sigma$ in $N$.  We define $w(s) = \la \dot\sigma(s),\bd_t\ra$ and write $\dot \sigma = w \bd_t + Z$ where $Z$ is tangent to the $t$-level set.  

\begin{prop}
We have 
$
1 = g_+(\dot \sigma,\dot\sigma) = w^2 + x^{-2} h_x(Z,Z). 
$
\end{prop} 

\begin{proof}
This is an immediate consequence of the fact that $\sigma$ is unit speed. 
\end{proof}

Write $\vert Z\vert^2 = \delta(Z,Z)$ and observe that $h_x(Z,Z) = \vert Z\vert^2 (1+O(x^3))$ and likewise $\vert Z\vert^2 = h_x(Z,Z) (1+O(x^3))$.  We will use this frequently in what follows.  Next, we show that if $\sigma$ is traveling toward the boundary it cannot ``turn around''. 

\begin{prop}
The function $s\mapsto t(\sigma(s))$ is convex. 
\end{prop}

\begin{proof}
We have 
\[
\frac{d^2}{ds^2} t(\sigma(s)) = \frac{d}{ds}w(s) =  \del^2_{g_+} t(\dot \sigma,\dot\sigma). 
\]
Write $\dot \sigma = w \bd_t + Z$ where $Z$ is tangent to the $t$-level set.  The Hessian satisfies 
\begin{gather*}
\del^2_{g_+}t(\bd_t,\bd_t) = 0,\\
\del^2_{g_+} t (Z,\bd_t) = Z(\vert \grad^{g_+} t\vert) = 0,\\
\del^2_{g_+} t(Z,Z) = \vert \grad^{g_+} t\vert A(Z,Z) = A(Z,Z)
\end{gather*} 
where $A$ is the second fundamental form of the $t$-level set. Since $g_+ = dt^2 + e^{2t}h_{e^{-t}}$, we have 
\[
A = \frac 1 2 \bd_t (e^{2t}h_{e^{-t}}) = e^{2t} h_{e^{-t}} + \frac 1 2 e^{2t} \bd_t h_{e^{-t}}.
\]
We have $x=e^{-t}$ and $h_x = \delta + O(x^3)$. 
We also have $\bd_t = -x\bd_x$ and hence 
\begin{align*}
A(Z,Z) &= x^{-2} \left( h_x(Z,Z) - \frac 1 2 x (\bd_x h_x)(Z,Z)\right) \\
&= x^{-2} \left(h_x(Z,Z) + O(x^3) \vert Z\vert^2\right) = x^{-2} h_x(Z,Z) (1+O(x^3)). 
\end{align*} 
Hence, using the unit speed identity, we deduce that 
\[
\del_{g+}^2 t(\dot\sigma,\dot\sigma) = A(Z,Z) = (1-w^2)(1+O(x^3)). 
\]
Thus, after potentially shrinking $\eps$, there is a uniform constant $c > 0$ such that 
\[
\del_{g+}^2 t(\dot\sigma,\dot\sigma) \ge c (1-w^2) 
\]
which proves the claimed convexity. 
\end{proof}

It is immediate from the previous proposition that $w(s) = \la \dot\sigma,\bd_t\ra$ is a non-decreasing function of $s$.  Moreover, since vertical lines are geodesics, if $w = 1$ or $w = -1$ at any point, then $w\equiv 1$ or $w\equiv -1$.  We now obtain a bound on the total possible horizontal movement of any geodesic contained in $N$. 

\begin{prop}
\label{horizontal-distance} 
Assume that $\sigma(s)$, $s\in [a,b]$ is a unit speed $g_+$ geodesic in $N$.  Write $\dot\sigma = w\bd_t + Z$ and let $\vert Z\vert$ be the Euclidean length of $Z$. Then we have 
\[
\int_a^b \vert Z(s)\vert\, ds \le C\eps
\]
where $C$ does not depend on $\sigma$. 
\end{prop}

\begin{proof}
If $\sigma$ is a vertical line this is obvious. So we can suppose $\sigma$ is not a vertical line, in which case we have $-1 < w < 1$ along $\sigma$. Then $w$ is a smooth, strictly increasing function of $s$ and so we can find a smooth inverse function $f$ such that $f(w(s)) = s$.  We then calculate the integral using the change of variables $s = f(w)$ which yields 
\begin{align*}
\int_a^b \vert Z(s)\vert\, ds = \int_{f^{-1}(a)}^{f^{-1}(b)} \vert Z(f(w))\vert f'(w)\, dw.
\end{align*} 
We have 
\[
f'(w) = \frac{1}{w'(s)} \le \frac{1}{c(1-w^2)}
\]
and by the unit speed identity we also have 
\[
\vert Z\vert = (\vert Z\vert^2)^{1/2} = \big(h_x(Z,Z) (1+O(x^3))\big)^{1/2} =  (1+O(x^3)) x \sqrt{1-w^2}\le 2\eps \sqrt{1-w^2}
\]
provided $\eps$ is small enough. 
Thus 
\[
\int_{f^{-1}(a)}^{f^{-1}(b)} \vert Z(f(w))\vert f'(w)\, dw \le \int_{-1}^1 \frac{2\eps }{c\sqrt{1-w^2}}\, dw = C\eps,
\]
and the result follows. 
\end{proof}

Now let $K' \subset K$ be a slightly smaller closed ball and then choose $\eps$ small enough that the number $C\eps$ from the previous proposition is smaller than the Euclidean distance from $\bd K'$ to $\bd K$.  Now consider the minimizing $g_+$ geodesic $\sigma$ from $p$ to a point $q\in N' = \{(x,z): x\in (0,\eps),\ z\in K'\}$.  By the previous proposition, the last entry of $\sigma$ into $N$ must occur through the top boundary $\{\eps\}\times K$ rather than the lateral boundary $(0,\eps)\times \bd K$. In particular, at the final point of entry this geodesic has $w\ge 0$ and from this point onward $w$ remains non-negative until $\sigma$ reaches $q$.  

\subsection{The Geodesic Equation} 
We now write out the geodesic equation. It is convenient to write the equation using cotangent vectors rather than tangent vectors, since this form of the equation involves only partial derivatives of the metric and not the Christoffel symbols. 
To that end, define $\xi = Z^\flat$ to be the cotangent vector identified with $Z$ via $g_+$.  Letting $z^1,z^2,z^3$ be Euclidean coordinates on $\R^3$, we have 
\[
Z=  Z^1 \frac{\bd}{\bd z^1} + Z^2 \frac{\bd}{\bd z^2} + Z^3 \frac{\bd}{\bd z^3}
\]
and hence the components of $\xi$ are given by 
\[
\xi_i = (g_+)_{ij} Z^j = x^{-2} (h_x)_{ij} Z^j
\]
where $i$ and $j$ run from $1$ to $3$.

Letting $t = z^0$ be a fourth coordinate, note that $\xi_i = \eta_i$ for $i=1,2,3$ where $\eta = \dot \sigma^\flat$. The geodesic equation therefore implies that the components of $\xi_i$ evolve according to the equation 
\[
\dot \xi_i = \frac 1 2 \bd_i (g_+)_{jk} Z^j Z^k = \frac 1 2 x^{-2} \sum_{j,k=1}^3   \bd_i (h_x)_{jk} Z^j Z^k 
\]
for $i=1,2,3$.  We have $\bd_i (h_x)_{jk} = O(x^3)$ and 
\begin{align*}
x^{-2} \sum_{j,k=1}^3 \vert Z^j\vert  \vert Z^k\vert   \le C x^{-2} \vert Z\vert^2 &= Cx^{-2} h_{x}(Z,Z) (1 + O(x^3)) \le C(1-w^2)
\end{align*} 
which implies that 
\[
\dot \xi_i = O(x^3(1-w^2)). 
\]
We write $\vert \xi\vert^2 = \delta^{ij} \xi_i \xi_j$. The above equation has the following consequence. 

\begin{prop}
At the endpoint $q$, we have $\vert \xi\vert \le C$ for a constant $C$ that is uniform over all $q\in N'$. 
\end{prop} 

\begin{proof}
Assume that $\sigma(a)$ is the last entry of $\sigma$ into $N$ and $\sigma(b) = q$.    We know that 
\[
h_x(Z,Z) = x^2(1-w^2)
\]
and so $\vert Z\vert \le C \eps \sqrt{1-w^2} \le C\eps$.  Moreover, we know from above that $x(\sigma(a)) = \eps$, and it follows that 
\[
\vert \xi_i(a) \vert \le \eps^{-2} \vert (h_x)_{ij}\vert \vert Z^j\vert \le \frac{C}{\eps}
\]
for $i=1,2,3$. Hence we obtain $\vert \xi\vert \le C/\eps$ at the point of entry. 

We have 
\[
\xi_i(b) = \xi_i(a) + \int_a^b \dot\xi_i(s) \, ds. 
\]
Writing $f(w(s)) = s$ and making the same change of variables $s = f(w)$ as in the proof of Proposition \ref{horizontal-distance}, we see that 
\[
\left\vert \int_a^b \dot\xi_i(s) \, ds\right\vert \le \int_{f^{-1}(a)}^{f^{-1}(b)} \vert \dot\xi_i(f(w))\vert f'(w) \, dw \le  \int_0^1 \frac{\vert \dot \xi_i(f(w))\vert }{c(1-w^2)}\, dw \le \frac{C}{c}. 
\]
Hence we have 
\[
\vert \xi (b)\vert \le \frac{C}{\eps} + \frac{C}{c} \le C(K)
\]
for a constant $C(K)$ that depends only on $K$. 
\end{proof} 

\subsection{Derivatives of the Distance Function} 

We now use the previous information about geodesics to study the derivatives of $r$.  Since $r$ is 1-Lipschitz, the gradient $\grad^{g_+} r$ exists at almost every point $q\in N'$. Choose such a point $q$ and let $\sigma$ be the minimizing unit speed $g_+$ geodesic from $p$ to $q$. As before, write $\dot \sigma = w \bd_t + Z$.  We have just seen that at the terminal point $q$ the cotangent vector $\xi = Z^\flat$ (determined by $g_+$) satisfies $\vert \xi\vert \le C$ where $C$ does not depend on $q$. This gives the following estimate for the derivatives of $r$. 

\begin{prop}
\label{horizontal-r-derivative}
At each point $q\in N'$ where $r$ is differentiable, we have $
\vert \bd_{z_i} r \vert \le C$. 
\end{prop}

\begin{proof}
We have 
\[
\bd_{z_i} r = g_+( \grad^{g_+} r, \bd_{z_i}) = g_+( \dot \sigma, \bd_{z_i}) = g_+( Z,\bd_{z_i}) = \xi(\bd_{z_i}). 
\]
Hence we have $\vert \bd_{z_i} r\vert = \vert \xi(\bd_{z_i})\vert \le \vert \xi\vert \vert \bd_{z_i}\vert \le C$. 
\end{proof} 

Finally, we record a simple identity for the $x$-derivative of $r$. 

\begin{prop}
At each point $q\in N'$ where $r$ is differentiable we have $\bd_x r = -w/x$. 
\end{prop}

\begin{proof}
We have 
\[
\bd_x r = g_+(\grad^{g_+} r,\bd_x) = g_+(\dot \sigma, -\bd_t/x) = -w/x,
\]
as needed. 
\end{proof} 

\subsection{Extending the Function to the Boundary} Recall that we have defined 
\[
u(x,z) = r(x,z) + \log x. 
\]
The function $u$ is locally Lipschitz on $N'$.  By Proposition \ref{horizontal-r-derivative}, we have $\vert \bd_{z_i} u\vert \le C$ almost-everywhere on $N'$. We also have 
\[
0 \le \bd_x u \le Cx
\]
almost-everywhere on $N'$ since 
\[
\bd_x u = \bd_x r + \frac 1 x = \frac{1-w}{x} 
\]
and 
\[
0 \le 1-w \le  1-w^2 = (x^{-2} h_x)(Z,Z) = (x^2 h_x^{-1})(\xi,\xi) \le C x^2 \vert \xi\vert^2 \le Cx^2. 
\]
These constants are allowed to depend on $K'$, but the same constant $C(K')$ works for all $(x,z) \in N'$. 

It now follows easily that $u$ extends to the boundary at infinity to give a Lipschitz function $u_0$ on $K'$.  Taking an exhaustion of $\R^3$ by compact sets and using the uniqueness of $x$, this defines a locally Lipschitz function $u_0$ on all of $\R^3$.  To summarize, we have the following. 

\begin{prop}\label{Prop: Lipschitz extension of u=r+logx}
There is a locally Lipschitz function $u_0$ on $\R^3$ such that we have a locally uniform expansion 
\[
u(x,z) = r(x,z) + \log x = u_0(z) + O(x^2).
\]
Moreover, the function $u$ satisfies $\vert \bd_{z_i} u\vert \le C$ and $0 \le \bd_x u \le Cx$ almost-everywhere in each collar $N = (0,\eps) \times K$ with $K\subset \R^3$ compact and $\eps$ sufficiently small for a constant $C = C(N)$. 
\end{prop}

\section{A Distributional Inequality} 
\label{distribution}

Recall that $r(q) = d_{g+}(q,p)$.  Since $\ric(g_+) = -3g_+$, it follows that $r$ satisfies 
\begin{equation}
\label{eq1} 
\lap_{g+} r \le 3\coth r
\end{equation} 
in the sense of distributions on $(\ins(X)-\{p\},g_+)$.  In particular, this distributional inequality holds even across the cut locus.  For a proof, see \cite[Corollary 2.12]{mantegazza2014distributional}.  

Define $F = e^{-r/2}$.  Then we have the following immediate consequence of this inequality.  Fix a compact set $K = \cl B(0,r)\subset \R^3$ and let $U$ be a neighborhood of $K$ on which there is a local geodesic defining function $x$. Choose a small $\eps > 0$ so that 
\[
N = \{(x,z): x\in (0,\eps),\ z\in K\}
\]
is contained in $U$. Note that $F < 1$ on $N$ as long as we choose $\eps$ small enough that $p\notin N$. 

\begin{prop}
The function $F$ satisfies 
\begin{equation}
\label{dist-F} 
\left(\lap_{g+} + \frac 5 4 \right)F  \ge -\frac{3F^5}{1-F^4}
\end{equation}
on $(N,g_+)$ in the sense of distributions. 
\end{prop}

\begin{proof}
Let $\eta\geq 0$ be a non-negative smooth function with compact support in $N$.  Let $(r_k)_{k\in \N}$ be a sequence of smooth approximations to $r$ with $r_k \to r$ uniformly and $\grad^{g_+} r_k \to \grad^{g_+} r$ in $L^1$ on the support of $\eta$.  Testing the inequality \eqref{eq1} against $e^{-r_k/2}\eta$ gives 
\[
\int r \lap_{g_+} (e^{-r_k/2} \eta)\, dv_{g_+} \le \int 3 \eta e^{-r_k/2}\coth r \, dv_{g_+}. 
\]
Since $r$ is Lipschitz, we can integrate by parts and use the fact that $\vert \grad^{g_+} r \vert = 1$ almost-everywhere to get 
\begin{align*}
& \int r \lap_{g_+} (e^{-r_k/2} \eta)\, dv_{g_+} \\
&\qquad = -\int e^{-r_k/2} \la \grad^{g_+} r, \grad^{g_+} \eta\ra \, dv_{g+} + \frac 1 2 \int \eta e^{-r_k/2} \la \grad^{g_+} r, \grad^{g_+} r_k\ra\, dv_{g_+}\\
&\qquad = -\int e^{-r_k/2} \la \grad^{g_+} r, \grad^{g_+} \eta\ra \, dv_{g+} + \frac 1 2 \int \eta e^{-r_k/2} \bigg[\la \grad^{g_+} r, \grad^{g_+} r_k - \grad^{g_+} r\ra + 1\bigg] \, dv_{g_+}.
\end{align*}
Hence we can send $k\to \infty$ to get 
\[
-\int e^{-r/2} \la \grad^{g_+} r, \grad^{g_+} \eta\ra \, dv_{g+} + \frac 1 2 \int \eta e^{-r/2} \, dv_{g_+} \le \int 3\eta e^{-r/2} \coth r\, dv_{g_+}. 
\]
Multiply by $-1/2$ and rearrange to get 
\begin{align*} 
\frac 1 2 \int e^{-r/2}\la \grad^{g_+} r, \grad^{g_+} \eta\ra \, dv_{g+} &\ge  \int \left( \frac 1 4 e^{-r/2}- \frac{3}{2} e^{-r/2}\coth r\right)\eta\, dv_{g_+}.
\end{align*} 
Finally, integrating by parts again shows that 
\[
\frac 1 2 \int e^{-r/2}\la \grad^{g_+} r, \grad^{g_+} \eta\ra \, dv_{g+} = \int e^{-r/2} \lap_{g_+} \eta\, dv_{g_+}. 
\]
Combining everything, we deduce that 
\begin{align*}
\int e^{-r/2} \lap_{g_+}\eta  + \frac 5 4 e^{-r/2} \eta \, dv_{g_+} &\ge \int \left( \frac 3 2 e^{-r/2}- \frac{3}{2} e^{-r/2}\coth r\right)\eta\, dv_{g_+} \\
&= -\int \frac{3 e^{-(5/2)r}\eta}{1 - e^{-2r}} \, dv_{g_+}. 
\end{align*}
This is precisely the distributional meaning of \eqref{dist-F}. 
\end{proof} 

Next we rewrite this distributional inequality in terms of the coordinate measure $dx\, dz$ on $N$. Define the function $\psi(x,z) = x^{-1/2} F = x^{-1/2} e^{-r/2}=e^{-u/2}$. Also define 
\[
J(x,z) = \sqrt{\det h_x}, \quad b(x,z) = \bd_x \log J.
\]
We compute that 
\[
\lap_{g_+} = x^2(\bd_x^2 + b\bd_x + \lap_{h_x}) - 2x\bd_x.
\]
Then we have the following equivalent form of \eqref{dist-F}. 

\begin{prop}
The function $\psi$ satisfies 
\begin{equation}
\label{eq5} 
\psi_{xx} - \frac{\psi_x}{x} + b\psi_x + \lap_{h_x} \psi + \frac{b}{2x} \psi \ge - \frac {3\psi^5}{1-x^2\psi^4}
\end{equation}
in the sense of distributions on $N$ with respect to $dx \, dz$. More precisely, this means that 
\[
\int \psi \left(\eta_{xx} +  \bd_x \left(\frac {\eta}{x}\right) - \bd_x(b\eta) + \frac{b}{2x} \eta + J \lap_{h_x}(J^{-1}\eta) \right) \, dx \, dz \ge - \int \frac{3\psi^5 \eta}{1-x^2\psi^4}\, dx\,dz
\]
for all smooth non-negative functions $\eta$ with compact support in $N$. 
\end{prop}

\begin{proof}
Fix any smooth function $\eta\ge 0$ compactly supported in $N$.  Note that the functions $x$ and $J$ are smooth.  Thus we can test \eqref{dist-F} against the function $\varphi = x^{3/2} J^{-1}\eta \ge 0$.  We calculate that 
\[
\varphi_x = \frac 3 2 x^{1/2} J^{-1} \eta - x^{3/2} J^{-1} b\eta + x^{3/2} J^{-1}\eta_x
\]
and 
\begin{align*} 
\varphi_{xx} &= \frac 3 4 x^{-1/2} J^{-1}\eta - 3 x^{1/2} J^{-1} b\eta + 3 x^{1/2} J^{-1} \eta_x + x^{3/2} J^{-1} b^2 \eta\\ &\qquad\quad   - x^{3/2} J^{-1} b_x \eta - 2 x^{3/2} J^{-1} b\eta_x + x^{3/2} J^{-1}\eta_{xx}. 
\end{align*} 
Using $dv_{g_+} = x^{-4} Jdxdz$, it then follows that 
\[
x^{1/2} \left(\lap_{g_+} + \frac 5 4\right) \varphi\, dv_{g_+} = (x^{-7/2} J) \left(x^2 \bd_x^2 + x^2 b \bd_x + x^2 \lap_{h_x} - 2x\bd_x + \frac 5 4\right) \varphi \, dx\,dz.
\]
Plugging in the above formulas and simplifying, we then get 
\[
x^{1/2} \left(\lap_{g_+} + \frac 5 4\right) \varphi\, dv_{g_+} = \left(\eta_{xx} + \bd_x \left(\frac \eta x \right) - \bd_x (b\eta) + \frac{b}{2x}\eta + J \lap_{h_x}(J^{-1}\eta)\right)\, dx \, dz.
\]
Thus from \eqref{dist-F}, we obtain 
\[
\int \psi \left(\eta_{xx} + \bd_x \left(\frac \eta x \right) - \bd_x (b\eta) + \frac{b}{2x}\eta + J \lap_{h_x}(J^{-1}\eta)\right)\, dx \, dz \ge - \int \frac{3 \psi^5 \eta}{1-x^2\psi^4}\, dx\, dz,
\]
as needed. 
\end{proof} 

Now we use this to prove a distributional inequality for $\psi_0:=e^{-u_0/2}$ on $\R^3$, where $u_0$ is given by Proposition \ref{Prop: Lipschitz extension of u=r+logx}. 

\begin{prop}
The function $\psi_0:=e^{-u_0/2}$ satisfies 
\[
-\lap \psi_0 \le 3\psi_0^5 
\]
on $\R^3$ in the sense of distributions. 
\end{prop}

\begin{proof}
Choose any smooth compactly supported function $\phi \ge 0$ on $\R^3$. We can assume that $K$ is chosen to contain the support of $\phi$.  Let $\eta(x) \ge 0$ be a smooth function with compact support in $(0,\eps)$ and then test \eqref{eq5} against the function $\phi(z)\eta(x)$. This gives 
\begin{align*} 
&\int_N \psi \left[\phi \eta_{xx} + \phi \bd_x \left(\frac \eta x\right) - \phi \bd_x (b\eta) + \eta J \lap_{h_x}(J^{-1}\phi) + \frac{b}{2x}\phi \eta\right]\, dx \, dz  \geq - \int_N \frac{3\psi^5 \phi \eta}{1-x^2\psi^4}\, dx\,dz. 
\end{align*} 
Define the functions 
\begin{gather*}
A_\phi(x) = \int_{\R^3} \psi(x,z)\phi(z)\, dz,\\
G_\phi(x) = \int_{\R^3} \frac{b}{2x}\psi \phi + b \phi \bd_x \psi + \psi J \lap_{h_x}(J^{-1}\phi) + \frac{3\psi^5 \phi}{1-x^2 \psi^4}\, dz.
\end{gather*}
The above equation can be re-written in the form 
\[
\int A_\phi(x) \left(\eta_{xx} + \bd_x \left(\frac \eta x\right)\right) + G_\phi(x)\eta\, dx \ge 0,
\]
and since $\eta$ was arbitrary it follows that 
\[
A_\phi'' - \frac{A_\phi'}{x} + G_\phi \ge 0
\]
in the sense of distributions on $(0,\eps)$. 

Next we study the behavior of $G_\phi$ as $x\to 0$.  Recall that $b= O(x^2)$ and $\psi = \psi_0 + O(x^2)$ and $\vert \psi_x\vert \le Cx$. Moreover, we have $h_x = \delta + O(x^3)$ with the analogous bounds on derivatives which implies that 
\[
J \lap_{h_x}(J^{-1}\phi) = \lap_\delta \phi + O(x^3). 
\]
Combining these estimates, it follows that 
\[
G_\phi(x) = \underbrace{\int \psi_0 \lap_\delta \phi + 3\psi_0^5 \phi \, dz}_{L_\phi} + q(x) 
\]
where $\vert q(x)\vert \le Cx$. 

To prove the proposition, we need to show that $L_\phi\ge 0$. Suppose to the contrary that $L_\phi < 0$. Then for a small $\delta > 0$, we have 
\[
A_\phi'' - \frac{A_\phi'}{x} \ge c > 0
\]
in the sense of distributions on $(0,\delta)$. Since $\vert \psi_x\vert \le Cx$, the function $A_\phi$ is locally Lipschitz and   
\[
B_\phi(x) = \frac{A_\phi'(x)}{x} 
\]
is a well-defined function in $L^\infty(0,\delta)$. Moreover, one can check that the above distributional inequality implies that 
\[
(B_\phi - c\log x)' \ge 0
\]
in the sense of distributions.  However, integrating this between Lebesgue points shows 
\[
B_\phi(x_2) - B_\phi(x_1) \ge c \log\left(\frac{x_2}{x_1}\right). 
\]
Letting $x_1\to 0$, the right hand side blows up while the left stays bounded and this is a contradiction.  Thus we have $L_\phi \ge 0$ as needed. 
\end{proof} 

\begin{rem}
    In the conformally compact case where $g_+$ is a metric on $B^4$ and $x^2 g_+$ induces the metric $h_0$ on $S^3$, the same argument applies, but the expansion of $b$ yields an extra term involving $R_{h_0}$. This  lets one get a distributional inequality for the conformal Laplacian $-8\lap_{h_0} + R_{h_0}$ on $S^3$. 
\end{rem}

\section{Volume Comparison} 
\label{volume}

Again select a closed ball $K$ in $\R^3$ and let $U$ be a neighborhood of $K$ in $X$ in which there is a smooth geodesic defining function $x$. Choose $\eps > 0$ for which 
\[
N = \{(x,z): x\in (0,\eps),\ z\in K\}
\]
is contained in $U$.  We can assume $\eps$ is small enough that $p\notin N$.  We have 
\[
\bd_x (r + \log x) \le C(K) x 
\]
on $N$ and thus (potentially shrinking $\eps$ further), we can suppose that 
\[
\bd_x r \le -\frac{1}{2x}
\]
on $N$. In particular, this implies that $r$ is a decreasing function of $x$ on $N$. Let $R_0 = \max_{z\in K} r(\eps,z)$.  Since $r\to \infty$ as $x\to 0$, it follows that for each $R > R_0$ and each $z\in K$, there will be a unique value $x_R(z)$ for which $r(x_R(z),z) = R$.  We now obtain an asymptotic expansion for $x_R(z)$. 

\begin{prop}
We have $x_R(z) = e^{u_0(z) - R}(1+O(e^{-2R}))$ as $R\to \infty$. 
\end{prop}

\begin{proof}
By definition, we have 
\[
r(x_R(z),z) = R, \quad u(x_R(z),z) = r(x_R(z),z) + \log x_R(z)
\]
and it follows that $u - R = \log x_R(z)$ which implies that 
\begin{equation}
\label{eq4} 
x_R(z) = e^{u(x_R(z),z) - R} = \left(e^{u_0(z) - R}\right)\left( e^{u(x_R(z),z)-u_0(z)}\right). 
\end{equation}
It remains to show that 
\begin{equation}
\label{eq3} 
 e^{u(x_R(z),z)-u_0(z)} = 1 + O(e^{-2R}). 
\end{equation}
We prove this in two stages.  First, recall that $u(x,z) - u_0(z) = O(x^2)$.  Since $0 < x_R(z) < \eps$ by construction, this implies that 
$
\vert u(x_R(z),z) - u_0(z)\vert \le C \eps^2
$
and therefore that 
\[
\vert e^{u(x_R(z),z) - u_0(z)} - 1 \vert \le C \eps^2. 
\]
Assuming $\eps$ is small enough, we can substitute this into \eqref{eq4} to obtain the improved bound 
\[
x_R(z) \le 2 e^{u_0(z) - R} \le C(K) e^{-R}. 
\] 
Finally this implies that 
\[
\vert u(x_R(z),z) - u_0(z)\vert \le C e^{-2R}
\]
and \eqref{eq3} follows. 
\end{proof} 

Next, consider the regions 
\[
\Omega_{R,K} = \{(x,z): z\in K,\, x_R(z) < x < \eps\}.
\]
Then by construction we have 
$
B_{g_+}(p,R) \cap N = \Omega_{R,K}. 
$
Next, we want to calculate the volume of these regions $\Omega_{R,K}$. 

\begin{prop}
\label{volume-est}
We have 
\[
\lim_{R\to \infty} e^{-3R} \vol_{g_+}(\Omega_{R,K}) = \frac 1 3 \int_{K} \psi_0^6\, dz. 
\]
\end{prop}

\begin{proof}
We have 
\[
dv_{g_+} = x^{-4} \sqrt{\det h_x} dxdz = x^{-4}(1+O(x^3)) dxdz
\]
and therefore 
\begin{align*}
\vol_{g_+}(\Omega_{R,K}) &= \int_{\Omega_{R,K}} x^{-4} + O(x^{-1}) \, dx\, dz \\
&= \int_K \left[\int_{x_R(z)}^\eps x^{-4} + O(x^{-1})\, dx\right]\, dz \\
&= \int_K \frac{1}{3x_R(z)^3}\, dz  -\int_K \frac{1}{3\eps^3}\, dz + E.
\end{align*}
Here the second term 
is a constant, 
and the error $E$ satisfies 
\begin{align*}
\vert E\vert &\le C \int_K \left[\int_{x_R(z)}^\eps x^{-1}\, dx\right]\, dz \\
&= C \int_K \bigg[\log(\eps) - \log(x_R(z))\bigg]\, dz\\
&\le C\left(1 + \int_K \vert \log x_R(z)\vert\, dz\right).
\end{align*} 
Since $x_R(z) \ge c e^{-R}$, it follows that 
\[
\vert E\vert \le C \left(1 + \int_K R\, dz\right) \le C R. 
\]
Finally, note that 
\begin{align*}
\int_K x_R(z)^{-3}\ dz &= \int_K e^{-3u_0 + 3R}(1+O(e^{-2R}))\, dz  
\end{align*} 
and so 
\begin{align*} 
e^{-3R} \int_K x_R(z)^{-3}\, dz &= \int_K e^{-3u_0} (1+O(e^{-2R}))\, dz \\
&= \int_K \psi_0^6 (1+O(e^{-2R}))\, dz \to \int_K \psi_0^6\, dz
\end{align*} 
as $R\to \infty$.  The result now follows. 
\end{proof}

This lets us calculate an upper bound on the $L^6$ norm of $\psi_0$. 

\begin{prop}
The function $\psi_0$ satisfies 
\[
\int_{\R^3}\psi_0^6\, dz \le \frac{\omega_3}{8}. 
\]
\end{prop}

\begin{proof}
Fix a compact $K\subset \R^3$. The previous proposition implies that 
\[
\frac{1}{3}\int_K \psi_0^6\, dz = \lim_{R\to\infty} e^{-3R} \vol_{g_+}\Omega_{R,K}. 
\]
On the other hand, we have $\Omega_{R,K} \subset B_{g_+}(p,R)$ and therefore 
\begin{align*}
\frac{1}{3}\int_K \psi_0^6 \, dz &\le \liminf_{R\to \infty} e^{-3R} \vol_{g_+} (B_{g_+}(p,R))\\
&\le \liminf_{R\to \infty} e^{-3R} \vol_{\mathbb H^4}(B(R)) = \frac{\omega_{3}}{24},
\end{align*}
where we used the Bishop-Gromov volume comparison to get from the first line to the second.  Since the upper bound does not depend on $K$, we can take an exhaustion of $\R^3$ by larger and larger compact sets and the result follows. 
\end{proof}

\begin{rem}
    This same argument applies to estimate the integral of $\psi_0^6$ in the conformally compact case on $B^4$. 
\end{rem}

Recall that a minimizer $f$ for the Euclidean $L^2$ Sobolev inequality satisfies $-\lap f = \lambda f^5$ where the constant $\lambda$ depends on the normalization of $f$. If we normalize so that 
\[
\int_{\mb R^3} f^6\, dz = \frac{\omega_3}{8}
\] 
then $\lambda = 3$.  We have just shown that $\psi_0$ satisfies corresponding inequalities $-\lap \psi_0 \le 3 \psi_0^5$ and 
\[
\int_{\R^3} \psi_0^6 \, dz \le \frac{\omega_3}{8}. 
\]
The Euclidean Sobolev inequality will then force these to be equalities. 

\begin{prop}
We have 
\[
\int_{\R^3} \psi_0^6\, dz = \frac{\omega_3}{8}. 
\]
\end{prop}

\begin{proof}
Choose test functions $\eta_R$ on $\R^3$ with $\eta_R \equiv 1$ on $B_R$ and $\eta_R \equiv 0$ outside $B_{2R}$ and $\vert \grad \eta_R\vert \le C/R$.  For fixed $R$, choose smooth functions $\psi_{0,k}$ such that $\psi_{0,k} \to \psi_0$ uniformly and $\grad \psi_{0,k} \to \grad \psi_0$ in $L^1$ on the support of $\eta_R$.  We can test the distributional inequality $-\lap \psi_0 \le 3\psi_0^5$ against $\psi_{0,k} \eta_R^2$ to get 
\[
-\int \psi_0 \lap (\psi_{0,k} \eta_R^2) \, dz \le 3 \int \psi_0^5 \psi_{0,k} \eta_R^2\, dz. 
\]
Then since $\psi_0$ is Lipschitz, we can integrate by parts to get 
\[
\int \eta_R^2 \la \grad \psi_0, \grad \psi_{0,k}\ra + 2 \psi_{0,k} \eta_R \la \grad \psi_0, \grad \eta_R\ra\, dz \le 3 \int \psi_0^5 \psi_{0,k}\eta_R^2\, dz. 
\]
Letting $k$ go to infinity, this implies that 
\[
\int \eta_R^2 \vert \grad \psi_0\vert^2 + 2 \psi_0 \eta_R \la \grad \psi_0, \grad \eta_R\ra \, dz \le 3 \int \psi_0^6 \eta_R^2 \, dz. 
\]
Completing the square, this implies that 
\[
\int \vert \grad(\eta_R \psi_0)\vert^2 \, dz \le 3\int \psi_0^6 \eta_R^2\, dz + \int \psi_0^2 \vert \grad \eta_R\vert^2\, dz. 
\]
Next, we apply the H{\"o}lder inequality to get 
\begin{align*}
\int \psi_0^2 \vert \grad \eta_R \vert^2\, dz &\le \left(\int_{B_{2R}-B_R} \psi_0^6\, dz\right)^{1/3} \left(\int_{B_{2R}-B_R} \vert \grad \eta_R\vert^3\, dz\right)^{2/3} \\
& \le C \left(\int_{B_{2R}-B_R} \psi_0^6\, dz\right)^{1/3} \to 0
\end{align*} 
as $R\to \infty$. We also have 
\[
\int_{B_R} \vert \grad \psi_0\vert^2 \, dz \le \int \vert \grad (\eta_R \psi_0)\vert^2\, dz 
\]
and therefore we can let $R\to \infty$ to deduce that 
\[
\int_{\R^3} \vert \grad \psi_0\vert^2\, dz \le 3\int_{\R^3} \psi_0^6\, dz. 
\]

 We can now apply the sharp Sobolev inequality to get 
\[
\frac{3 \omega_3^{2/3}}{4} \left(\int_{\R^3} \psi_0^6\, dz \right)^{1/3}  \le  \int_{\R^3} \vert \grad \psi_0\vert^2\, dz.
\]
Abbreviating 
\[
M = \int_{\R^3} \psi_0^6\, dz, 
\]
we can combine this with the previous bound to get 
\[
\frac{3 \omega_3^{2/3}}{4} M^{1/3} \le 3M 
\]
which implies that 
\[
M \ge \frac{\omega_3}{8}. 
\]
But we have already demonstrated that 
\[
M \le \frac{\omega_3}{8}
\]
and the equality follows. 
\end{proof}

Finally, returning to the volume relationship, we can now prove equality in the Bishop-Gromov volume comparison, which then implies the isometry to $\mathbb H^4$. 

\begin{prop}\label{Prop: rigidity}
The manifold $(X,g_+)$ is isometric to $(\mathbb H^4,g_+)$. 
\end{prop}

\begin{proof}
Recall that we have already shown that 
\begin{align*}
\frac 1 3 \int_{\R^3} \psi_0^6\, dz &\le \liminf_{R\to \infty} e^{-3R} \vol_{g_+}( B_{g_+}(p,R))\\
& = \liminf_{R\to \infty} \left[e^{-3R}\vol_{\mathbb H^4} (B(R)) \frac{\vol_{g_+}( B_{g_+}(p,R))}{ \vol_{\mathbb H^4} (B(R))}\right] \\
&= \frac{\omega_3}{24} \liminf_{R\to \infty} \frac{\vol_{g_+}( B_{g_+}(p,R))}{ \vol_{\mathbb H^4} (B(R))} \le \frac{\omega_3}{24}. 
\end{align*} 
On the other hand, we have just demonstrated that 
\[
\int_{\mb R^3} \psi_0^6\, dz = \frac{\omega_3}{8}
\]
and it follows that all of the above inequalities must be equalities.  In particular, we have 
\[
\liminf_{R\to \infty} \frac{\vol_{g_+}( B_{g_+}(p,R))}{ \vol_{\mathbb H^4} (B(R))} = 1. 
\]
It is well-known that this implies $(X,g_+)$ is isometric to $(\mathbb H^4,g_+)$. 
\end{proof}

\part{The existence of equivariant Poincar{\'e}-Einstein fillings}
\section{Preliminaries for equivariant Poincar{\'e}-Einstein fillings}

Given any $p\in\{2,3,\dots\}$ and coprime $1\leq q<p$, let
\[ \mk r_{p,q}: (z_1,z_2) \mapsto (e^{2\pi i/p}\cdot z_1, e^{2\pi iq/p}\cdot z_2) \]
be a diffeomorphism on the unit $3$-sphere $S^3=\{(z_1,z_2)\in\mb C^2 : |z_1|^2+|z_2|^2=1\}$. 
Then 
\[ \mb Z_p:=\langle \mk r_{p,q}\rangle \]
is a free isometric action on $S^3$ so that $S^3/\mb Z_p = L(p,q)$. 
In this paper, we use the following notations:
\begin{itemize}
    \item $B^4$, $\overline{B}^4$: the open/closed unit Euclidean ball in $\mb R^4$;
    \item $Met(S^3)$: the space of smooth metrics on $S^3$;
    \item $Met_{\mb Z_p}(S^3)$: the space of $\mb Z_p$-invariant smooth metrics on $S^3$;
    \item $\mc C$: the space of conformal classes of metrics in $Met(S^3)$;
    \item $h_0$: the round metric on the unit sphere $S^3=\bd \overline{B}^4$. 
\end{itemize}

\subsection{Equivariant parameterization of $\mc C$}

\begin{lem}\label{Lem: action on conformal classes}
    The diffeomorphism $\mk r_{p,q}$ on $S^3$ induces a well-defined smooth map 
    \[\mk r_{p,q}:  [h]\mapsto [\mk r_{p,q}^*h]\] 
    on $\mc C$. 
    Let $\mb Z_p:=\langle\mk r_{p,q}\rangle$ be the action on $\mc C$ generated by the above map $\mk r_{p,q}$, and denote by 
    \[\mc C_{\mb Z_p}:=\{[h]\in \mc C: [\mk r_{p,q}^*h]=[h]\}\] 
    the $\mb Z_p$-fixed-point set. 
    Then, $[h]\in \mc C_{\mb Z_p}$ if and only if there is a $\mb Z_p$-invariant metric $\bar h\in [h]$. 
\end{lem}
\begin{proof}
    Given any $h'\in [h]$, there exists a function $f>0$ on $S^3$ so that $h'=f^2h$. 
    Hence, $\mk r_{p,q}^*h'=(f\circ \mk r_{p,q})^2 \cdot \mk r_{p,q}^*h\in [\mk r_{p,q}^*h]$, which implies that the map $\mk r_{p,q}:  [h]\mapsto [\mk r_{p,q}^*h]$ is well-defined on $\mc C$. 
    Clearly, if $[h]\in\mc C$ has a $\mb Z_p$-invariant representative $\bar h\in [h]$, then $[h]=[\bar h]=\mk r_{p,q}\cdot [\bar h]\in \mc C_{\mb Z_p}$. 
    Meanwhile, for any $[h]\in \mc C_{\mb Z_p}$ and $1\leq k\leq p$, we have $(\mk r_{p,q}^k)^* h = f_k^2\cdot h $ for some function $f_k>0$. 
    Then, 
    \begin{align}
        \bar h:= \frac{1}{p}\sum_{k=1}^p (\mk r_{p,q}^k)^*h = \frac{\sum_{k=1}^p f_k^2}{p}\cdot h 
    \end{align}
    is a $\mb Z_p$-invariant metric as a representative of $[h]$.
\end{proof}

Given any $h\in Met(S^3)$, the quotient $d\vol_{h_0}/d\vol_h$ of the volume densities is a smooth positive function on $S^3$. 
Hence, every $[h]\in\mc C$ admits a unique {\em unit-density representative} $\hat h$ given by 
\begin{align}\label{Eq: unit-density representative}
    \hat h:= \left(\frac{d\vol_{h_0}}{d\vol_h}\right)^{\frac{2}{3}} \cdot h 
\end{align}
so that $d\vol_{h_0}=d\vol_{\hat h}$. 
In particular, if $[h]\in\mc C_{\mb Z_p}$, then the unit-density representative $\hat h$ is a $\mb Z_p$-invariant metric on $S^3$. 

Using the unit-density representatives, we can identify the space $\mc C$ with the space of unit-density smooth metrics on $S^3$.  
Moreover, we show in the following lemma that $\mc C$ can be parameterized by a linear space $\mc A$ so that the pull-back action $[h]\mapsto [\mk r_{p,q}^*h] $ on $\mc C$ is represented by a linear action on the parameter space. 

\begin{lem}\label{Lem: parameterization of boundary conformal classes}
    Let $\mc A$ be the space of $h_0$-self-adjoint trace-free smooth endomorphism fields on $S^3$. 
    Namely, for any $A\in \mc A$, $A=\sqcup_{x\in S^3}A_x$ is a smooth family of $h_0$-self-adjoint trace-free linear maps with $A_x\in\text{End}(T_xS^3)$. 
    Then there exists a diffeomorphism 
    \[F: \mc C\to \mc A,\qquad F([h]):= \log(h_0^{-1}\hat h)=\log\left(\left((h_0)^{ik}(\hat h)_{kj}\right)_{3\times 3}\right),\] 
    so that for any $[h]\in \mc C$ with the unit-density representative $\hat h$, $A:=F([h])\in \mc A$ satisfies
    \[ \hat h(u,v)= h_0(e^{A_x}\cdot u,v)=h_0(u,e^{A_x}\cdot v),\qquad \mbox{$\forall u,v\in T_xS^3$ and $x\in S^3$}, \]
    where $A_x:=A|_x: T_xS^3\to T_x S^3$, and $e^{A_x}:=\Id_{T_xS^3}+\sum_{k=1}^\infty\frac{A_x^k}{k!}$. 
    Additionally, 
    \[  F( [\mk r_{p,q}^*h])|_x = \left[(d\mk r_{p,q})^{-1}\circ F([h])|_{\mk r_{p,q}(x)}\circ d\mk r_{p,q} \right]_x  ,\qquad\forall x\in S^3.\]
    In particular, $F(\mc C_{\mb Z_p})$ is a linear subspace of $\mc A$. 
\end{lem}

\begin{rem}\label{Rem: equivariant parameterization of boundary conformal classes}
    Let $\mk r_{p,q}$ act on $\mc A$ linearly by $\mk r_{p,q}\cdot A:=(d\mk r_{p,q})^{-1}\circ A\circ d\mk r_{p,q}$ for any $A\in\mc A$. 
    Then, the above lemma shows that the diffeomorphism $F:\mc C\to \mc A$ is $\mb Z_p=\langle \mk r_{p,q}\rangle$ equivariant. 
    For simplicity, we denote by $\mc A_{\mb Z_p}$ the $\langle \mk r_{p,q}\rangle$-fixed-point set of $\mc A$. 
\end{rem}


\begin{proof}
    For any $x\in S^3$, we can choose an $h_0$-orthonormal basis $\{e_1,e_2,e_3\}$ of $T_xS^3$ so that $\hat h_x:=\hat h|_x$ is represented by a diagonal matrix $\diag(\lambda_1,\lambda_2,\lambda_3)$ with $\lambda_1,\lambda_2,\lambda_3>0$. 
    Since $d\vol_{h_0}=d\vol_{\hat h}$, we conclude that $\lambda_1\cdot \lambda_2\cdot \lambda_3=1$.  Hence, the diagonal matrix $\diag(\ln(\lambda_1), \ln(\lambda_2), \ln(\lambda_3))$ determines an $h_0$-self-adjoint trace-free endomorphism $A_x$ of $T_xS^3$. 
    For any $u=\sum_{i=1}^3 u_ie_i,v=\sum_{i=1}^3v_ie_i\in T_xS^3$, one easily sees that 
    \[\hat h_x(u,v) = (u_1,u_2,u_3)\cdot e^{\diag(\ln(\lambda_1), \ln(\lambda_2), \ln(\lambda_3))} \cdot (v_1,v_2,v_3)^T= (h_0)_x(e^{A_x}u,v) = (h_0)_x(u,e^{A_x}v). \]
    Hence, every $[h]\in \mc C$ uniquely determines an element $A=F([h])\in \mc A$ so that $F$ is smooth. 
    Conversely, every $A\in \mc A$ also uniquely and smoothly determines a unit-density metric $h:=h_0(e^A\cdot,\cdot)$. 
    This gives the diffeomorphism $F:\mc C\to\mc A$. 

    For any $x\in S^3$ and $u,v\in T_xS^3$, a direct computation shows that 
    \begin{align*}
        ((\mk r_{p,q})^*\hat h)_x (u,v) &= \hat h_{\mk r_{p,q}(x)} \left( d\mk r_{p,q}(u), d\mk r_{p,q}(v) \right) = (h_0)_{\mk r_{p,q}(x)} \left( e^{A_{\mk r_{p,q}(x)}}\circ d\mk r_{p,q} (u), d\mk r_{p,q}(v) \right)
        \\
        &= (\mk r_{p,q}^* h_0)_x \left( (d\mk r_{p,q})^{-1}\circ e^{A_{\mk r_{p,q}(x)}}\circ d\mk r_{p,q} (u), v \right)
        \\
        &= (h_0)_x \left( e^{(d\mk r_{p,q})^{-1}\circ A_{\mk r_{p,q}(x)}\circ d\mk r_{p,q}} (u), v \right),
    \end{align*}
    where we used the facts that $\mk r_{p,q}^*h_0=h_0$ and $e^{(d\mk r_{p,q})^{-1}\circ A_{\mk r_{p,q}(x)}\circ d\mk r_{p,q}}=(d\mk r_{p,q})^{-1}\circ e^{A_{\mk r_{p,q}(x)}}\circ d\mk r_{p,q}$ in the last equality. \
    Hence, we conclude that $F([\mk r_{p,q}^* h] )|_x= [(d\mk r_{p,q})^{-1}\circ F([h])|_{\mk r_{p,q}(x)}\circ d\mk r_{p,q} ]_x$. 
    This shows that the $\langle\mk r_{p,q}\rangle$-action on $\mc C$ induces a linear $\mb Z_p$-action on $\mc A$ given by 
    \[\mk r_{p,q}: A\mapsto (d\mk r_{p,q})^{-1}\circ A\circ d\mk r_{p,q},\qquad \forall A \in \mc A, \] 
    so that $F$ is a $\mb Z_p$-equivariant diffeomorphism. 
    Therefore, $F(\mc C_{\mb Z_p})$ is the $\mb Z_p$-fixed-point set of $\mc A$, which must be a linear subspace. 
\end{proof}

\subsection{Poincar{\' e}-Einstein fillings}

\begin{definition}
    A smooth function $\rho\in C^\infty(\overline{B}^4)$ is called a {\em defining function} on $\overline{B}^4$ if 
    \[\rho^{-1}(0)=\bd \overline{B}^4\qquad {\rm and}\qquad d\rho\neq 0  \quad\mbox{on $\bd \overline{B}^4$}.\] 
    Let $g_+$ be a complete metric on $B^4=\interior(\overline{B}^4)$. 
    If there exists a defining function $\rho$ on $\overline{B}^4$ so that the conformal metric 
    \[ \bar g:=\rho^2 g_+ \]
    extends to a $C^{m,\alpha}$ Riemannian metric on $\overline{B}^4$, then we say $g_+$ is {\em $C^{m,\alpha}$ conformally compact} on $B^4$. 
    Additionally, $\bar g$ is said to be the {\em $C^{m,\alpha}$ conformal compactification} of $g_+$ associated to the defining function $\rho$, and the induced metric $h=\bar g|_{\bd \overline{B}^4}$ is said to be the {\em boundary metric} associated to the compactification $\bar g$. 
\end{definition}

Note that defining functions are not unique. 
Hence, there are many conformal compactifications of a given metric $g_+$. 
Nevertheless, the conformal class $[h]$ of the boundary metric is uniquely determined by $(B^4, g_+)$.
In particular, $[h]$ is called the {\em conformal infinity} of $g_+$.

\begin{definition}
    Suppose the metric $g_+$ is $C^{m,\alpha}$ conformally compact on $B^4$ with $m\geq 2$. 
    If 
    \[ \Ric_{g_+} = -3 g_+ ,\]
    then $g_+$ is called a {\em conformally compact Einstein metric on $B^4$} (of $C^{m,\alpha}$ regularity), and $(B^4, g_+)$ is called a {\em conformally compact Einstein manifold} or a {\em Poincar{\'e}-Einstein manifold}. 
\end{definition}

Given a conformally compact Einstein metric $g_+$ on $B^4$ with a defining function $\rho>0$, we see from the formula of the conformal Ricci curvature that  
\[ -3\bar g=-3\rho^2g_+=\rho^2\Ric_{g_+}=\rho^2\Ric_{\bar g} +2\rho\nabla^2_{\bar g}\rho + (\rho\Delta_{\bar g} \rho - 3|d\rho|^2)\bar g ,\]
where $\bar g=\rho^2g_+$. 
On the boundary $\bd \overline{B}^4$, we have $\rho=0$, and thus $|d\rho|_{\bar g}=1$ on $\bd \overline{B}^4$. 
If the defining function $\rho$ also satisfies that 
\[\mbox{$|d\rho|_{\bar g}=1$ in a neighborhood of $\bd \overline{B}^4$},\] 
then $\rho$ is said to be a {\em geodesic defining function}.

Additionally, we also use the following notations:
\begin{itemize}
    \item $E_{PE}$: the space of $C^\infty$ conformally compact Einstein metrics $g_+$ on $B^4$ endowed with the $C^\infty$ topology induced from the map $g_+\mapsto \rho^2g_+$ into $C^\infty(\overline{B}^4, S^2T^*\overline{B}^4)$;
    \item $\Diff(\overline{B}^4)$: the group of $C^\infty$ diffeomorphisms $\Phi$ of $\overline{B}^4$;
    \item $\Diff_{\bd}(\overline{B}^4)$: the group of $C^\infty$ diffeomorphisms $\Phi$ of $\overline{B}^4$ with $\Phi|_{\bd \overline{B}^4}=\Id_{\bd \overline{B}^4}$; 
    \item $\mc E:= E_{PE}/\Diff_{\bd}(\overline{B}^4)$: the space of equivalence classes $[g_+]_{PE}$ of $C^\infty$ conformally compact Einstein metrics on $B^4$ endowed with the $C^\infty$ quotient topology, 
\end{itemize}
where $\Diff_{\bd}(\overline{B}^4)$ acts on $E_{PE}$ naturally by pullback. 
Moreover, there is a natural boundary map
\begin{align}\label{Eq: boundary map}
    \Pi: \mc E\to \mc C, \qquad \Pi([g_+]_{PE}) := [h],
\end{align}
where $[h]$ is the conformal class of the boundary metric associated to the conformal compactification of $g_+$.

Note that the hyperbolic metric on $B^4$ is a conformally compact Einstein metric. 
Hence, $\mc E$ is a non-empty set. 

\begin{proposition}\label{prop: action on the filling space}
    Using the above notations, set 
    \[ \wha \Pi := F\circ \Pi ~:~ \mc E\to \mc A ,\]
    and extend the diffeomorphism $\mk r_{p,q}$ from $S^3$ to $\overline{B}^4$ by 
    \[\mk r_{p,q}(z_1,z_2):= (e^{2\pi i/p}\cdot z_1, e^{2\pi iq/p}\cdot z_2),\qquad\forall (z_1,z_2)\in \mb C^2~\text{with}~|z_1|^2+|z_2|^2\leq 1,\]
    which further induces a well-defined continuous map on $\mc E$ given by
    \[\mk r_{p,q}: [g_+]_{PE}\mapsto [\mk r_{p,q}^*g_+]_{PE}.\] 
    Let $\mb Z_p:=\langle \mk r_{p,q} \rangle$ be the action on $\mc E$ generated by the above map $\mk r_{p,q}$, and denote by 
    \[ \mc E_{\mb Z_p} :=\left\{[g_+]_{PE}\in \mc E: [\mk r_{p,q}^*g_+]_{PE}=[g_+]_{PE}\right\} \]
    the $\mb Z_p$-fixed-point set of $\mc E$. 
    Then, 
    \begin{itemize}
        \item[(i)] $\Pi$ and $\wha \Pi$ are $\mb Z_p$-equivariant maps;
        \item[(ii)] for any $[g_+]_{PE}\in \mc E_{\mb Z_p}$, there exists a $C^\infty$ diffeomorphism $\phi$ on $\overline{B}^4$ with $\phi|_{\bd\overline{B}^4}=\mk r_{p,q}$ so that $\phi^*g_+=g_+$, and $\phi$ generates a $C^\infty$ $\mb Z_p$-action on $\overline{B}^4$. 
    \end{itemize}
\end{proposition}


\begin{proof}
    Suppose $g_+'\in [g_+]_{PE}$. 
    Then there exists $\Phi\in \Diff_\bd(\overline{B}^4)$ so that $g_+'=\Phi^*g_+$. 
    For the $C^\infty$ diffeomorphism $\mk r_{p,q}$ on $\overline{B}^4$, we have 
    \[\mk r_{p,q}^* g_+' =\mk r_{p,q}^*\Phi^*g_+ = (\mk r_{p,q}^{-1}\circ\Phi\circ\mk r_{p,q})^* (\mk r_{p,q}^* g_+).  \]
    Combining with the fact that $(\mk r_{p,q}^{-1}\circ\Phi\circ\mk r_{p,q})|_{\bd\overline{B}^4}=\Id_{\bd \overline{B}^4}$, we conclude that $\mk r_{p,q}^*g_+'\in [\mk r_{p,q}^*g_+]_{PE}$, which implies that the map $\mk r_{p,q}: [g_+]_{PE}\mapsto [\mk r_{p,q}^*g_+]_{PE}$ on $\mc E$ is well-defined. 

    Consider the $\langle\mk r_{p,q}\rangle$-actions on $\mc E$ and $\mc C$. 
    Given $[g_+]_{PE}\in \mc E$, suppose $\rho$ is a defining function of $g_+$ and $\bar g:=\rho^2g_+$ is the conformal compactification of $g_+$. 
    Then, $\mk r_{p,q}^*\bar g= (\rho\circ\mk r_{p,q})^2\cdot( \mk r_{p,q}^*g_+)$ is the conformal compactification of $\mk r_{p,q}^*g_+$ with the defining function $\rho\circ\mk r_{p,q}$. 
    Thus, $\Pi(\mk r_{p,q}\cdot[g_+]_{PE})=\Pi([\mk r_{p,q}^* g_+]_{PE})=[\mk r_{p,q}^*\bar g|_{\bd \overline{B}^4}]=\mk r_{p,q}\cdot [\bar g|_{\bd \overline{B}^4}]=\mk r_{p,q}\cdot \Pi([g_+]_{PE})$. 
    This together with Remark \ref{Rem: equivariant parameterization of boundary conformal classes} shows that $\Pi$ and $\wha \Pi$ are $\mb Z_p$-equivariant maps. 

    If $[g_+]_{PE}\in \mc E_{\mb Z_p}$, then there exists $\Phi\in \Diff_\bd(\overline{B}^4)$ so that $\mk r_{p,q}^* g_+=\Phi^*g_+$. 
    In particular, $\phi:= \mk r_{p,q}\circ \Phi^{-1}$ is a $C^\infty$ diffeomorphism on $\overline{B}^4$ with $\phi|_{\bd\overline{B}^4}=\mk r_{p,q}$ so that $\phi^*g_+=g_+$. 
    
    Consider the map $J:=\phi^p$ defined by compounding $\phi$ with itself for $p$ times. 
    Then, 
    \[J|_{\bd \overline{B}^4}=\mk r_{p,q}^p=\Id_{\bd\overline{B}^4} \qquad\text{and}\qquad J^*g_+=g_+.\] 
    We claim that $J=\Id_{\overline{B}^4}$, and thus $
    \phi$ generates a $C^\infty$ $\mb Z_p$-action on $\overline{B}^4$. 
    
    Indeed, let $\hat h$ be the unit-density representative of $[h]:=\Pi([g_+]_{PE})$. 
    Then, $\hat h$ is a $\langle \mk r_{p,q}\rangle$-invariant metric, and $\hat h$ determines a geodesic defining function $\rho$ with $\hat h= (\rho^2 g_+)|_{\bd\overline{B}^4}$ so that $\rho$ is {\em unique} in a neighborhood of $\bd\overline{B}^4$. 
    Define $\rho':=\rho\circ J$. 
    Since $J$ is a $g_+$-isometry, $\rho'^2\cdot g_+=J^*(\rho^2\cdot g_+)= J^*\bar g$, where $\bar g:=\rho^2g_+$ is the conformal compactification of $g_+$ associated to $\rho$. 
    Therefore, $|d\rho'|_{\rho'^2g_+}=|J^*d\rho|_{J^*\bar g}=|d\rho|_{\bar g}=1$ in a neighborhood of $\bd\overline{B}^4$ since $\rho$ is a geodesic defining function. 
    Additionally, since $J|_{\bd \overline{B}^4}=\Id_{\bd \overline{B}^4}$, the boundary metric of $g_+$ associated to $\rho'$ is $(\rho'^2\cdot g_+)|_{\bd \overline{B}^4}=(J^*\bar g)|_{\bd \overline{B}^4}=(J|_{\bd \overline{B}^4})^* \hat h=\hat h$. 
    In particular, we see that $\rho'$ is another geodesic defining function associated to the same $g_+$ and $\hat h$. The uniqueness shows that $\rho=\rho'$ and $J^* \bar g=\bar g$ in a neighborhood of $\bd \overline{B}^4$. 
    Next, in a small neighborhood of $\bd \overline{B}^4$, since $J$ preserves $\bar g$ and $\rho$, $J$ also preserves the flow generated by $\nabla_{\bar g}\rho$. 
    Recall that $J|_{\bd\overline{B}^4}=\Id_{\bd\overline{B}^4}$. 
    Hence, $J$ can only be the identity map in a neighborhood of $\bd\overline{B}^4$ to preserve the flow of $\nabla_{\bar g}\rho$. 
    Finally, take a small ball $B_r(x)\subset B^4$ so that $J=\Id$ in $B_r(x)$. 
    The completeness of $g_+$ on $B^4$ implies that for any $y\in B^4$, there exists $v\in T_xB^4$ with $y=\exp_x(v)$, where $\exp$ denotes the exponential map under the metric $g_+$. 
    Since $J$ is a $g_+$-isometry and $dJ|_{T_xB^4}=\Id_{T_xB^4}$, we conclude that $J(y)=\exp_{J(x)}(dJ(v))=y$. 
    This shows $J=\Id_{\overline{B}^4}$. 
\end{proof}

\section{Compactness theorem}

In this section, we show a compactness result based on the following compactness theory of Chang-Ge \cite{chang2025problem}. 
We refer to \cite{changGe2018compactness}\cite{changGeQing2020compactness} for related compactness results. 

\begin{theorem}[Chang-Ge \cite{chang2025problem}*{Theorem I}]\label{Thm: compactness of chang-ge}
    Fix $m\geq 6$, $\alpha\in (0,1)$, and an oriented $4$-manifold $X$ with $H_2(X;\mb Z)=H^2(X;\mb Z)=0$. 
    Let $\{g^i_+\}$ be a family of conformally compact Einstein metrics on $X$ with conformal infinity $(\bd X, [h_i])$. 
    Suppose 
    \begin{itemize}
        \item[(i)] $Y(\bd X, [h_i])\geq C_0>0$ for some uniform constant $C_0>0$;
        \item[(ii)] the boundary metrics $h_i$ have positive scalar curvature, and $\{h_i\}$ is compact in the $C^{m,\alpha}$ norm;
        \item[(iii)] $\int_{X}|W_{g_+^i}|^2 \leq C_1$ for some uniform constant $C_1>0$,
    \end{itemize}
    where $W_{g_+^i}$ is the Weyl curvature tensor of the metric $g_+^i$. 
    Then, the family of adapted scalar flat compactifications $\{\bar g_i\}$ is precompact in the $C^{m,\mu}$ norm for any $\mu\in (0,\alpha)$ up to modulo diffeomorphisms fixing $\bd X$. 
\end{theorem}

The main result of this section is the following theorem.
\begin{theorem}\label{Thm: compactness}
    Take $X=B^4$, $\bd X=S^3$, $m\geq 6$, and $\alpha\in (0,1)$. 
    Suppose $\{g^i_+\}$ is a family of conformally compact Einstein metrics on $X$ with conformal infinity $(\bd X, [h_i])$ so that (i) and (ii) in Theorem \ref{Thm: compactness of chang-ge} are satisfied. 
    Then there exists a constant $C_1>0$ so that every adapted scalar flat compactification $\bar g_i$ of $g_+^i$ satisfies 
    \[\sup_{X}\left( |\Rm_{\bar g_i}| + |\nabla \Rm_{\bar g_i}|^{2/3} \right) \leq C_1 \qquad\text{and}\qquad \qquad\vol(X,\bar g_i)\leq C_1.  \]
    In particular, $\int_{X} |W_{g_+^i}|^2=\int_{X} |W_{\bar g_i}|^2$ is uniformly bounded from above, and $\{\bar g_i\}$ is precompact in the $C^{m,\mu}$ norm for any $\mu\in (0,\alpha)$ up to modulo diffeomorphisms fixing $\bd X$ by Theorem \ref{Thm: compactness of chang-ge}. 
\end{theorem}

\subsection{The scalar flat adapted compactification}

We begin with some basic facts about the scalar-flat adapted compactification. 

Let $(B^4, g_+)$ be a conformally compact Einstein manifold with conformal infinity $[h]$ at $\bd \overline{B}^4=S^3$. 
Suppose the boundary metric $h$ has positive scalar curvature. 
By \cite{chang2025problem}*{Lemma 2.6}, there exists a defining function $\rho$ so that $\bar g:=\rho^2g_+$ extends to a $C^{m,\alpha}$ metric on $\overline{B}^4$ satisfying 
\begin{align}\label{Eq: adapted scalar-flat compactification}
    \bar g|_{\bd\overline{B}^4}=h,\qquad R_{\bar{g}}\equiv 0,\qquad |\nabla \rho|_{\bar g}\leq 1 ~\text{in}~B^4,\qquad\text{and}\quad H_{\bar g}>0~\text{on}~\bd\overline{B}^4, 
\end{align}
where $H_{\bar g}$ is the mean curvature on the boundary with respect to the outward unit normal. 
The metric $\bar g$ is called the {\em scalar flat adapted compactification} associated with $h$. 

\begin{lemma}\label{Lem: computation for scalar flat adapted metric}
    Given a scalar flat adapted compactification $\bar g=\rho^2 g_+$, consider the quantity 
    \begin{align}\label{Eq: P quantity}
        P:=\frac{1-|d\rho|_{\bar g}^2}{\rho}\geq 0.
    \end{align}
    Then, $P$ extends to a $C^{m-1,\alpha}$ function on $\overline{B}^4$, and
    \begin{align}
        \nabla^2\rho &= -\frac{\rho}{2}\Ric_{\bar g} - \frac{P}{2}\bar g, \label{Eq: scalar falt adapted metric - 3}
        \\
        \Delta\rho &= -2P, \label{Eq: scalar falt adapted metric - 4}
        \\
        \nabla P &= \Ric_{\bar g}(\nabla \rho, \cdot). \label{Eq: scalar falt adapted metric - 5}
    \end{align}
    Additionally, the boundary $(\bd\overline{B}^4,h)$ is umbilic with $H_{\bar g}=\frac{3}{2}P$ and 
    \begin{align}\label{Eq: scalar falt adapted metric - 6}
        R_h = -2\Ric_{\bar g}(\nu,\nu) +\frac{2}{3} H_{\bar g}^2,
    \end{align}
    where $\nu$ is the outward pointing unit normal of $\bd\overline{B}^4$. 
\end{lemma}

\begin{proof}
 By the conformal change formula of Ricci curvature in dimension $4$, we have 
 \[
 \ric_{g_+}=\ric_{\bar{g}}+2\rho^{-1}\nabla^2\rho+\left(\rho^{-1}\Delta\rho-3\rho^{-2}|\nabla\rho|^2\right)\bar{g}.
 \]

 Since $R_{\bar{g}}=0$ and $g_+$ is PE, we have 
 \[
 \Delta \rho=2\rho^{-1}\left(|\nabla \rho|^2-1\right),
 \]
 which is equivalent to \eqref{Eq: scalar falt adapted metric - 4}. This also proves the smoothness of $P$.

Combining \eqref{Eq: scalar falt adapted metric - 4} with the conformal change formula, we have 
\[
\nabla^2\rho
=-\frac{\rho}{2}\ric_{\bar{g}}-\frac{P}{2}\bar{g},
\]
which is \eqref{Eq: scalar falt adapted metric - 3}.

Differentiating $P$, we have 
\begin{align*}
    \nabla P=&\frac{-2\nabla^2\rho(\nabla \rho,\cdot)\rho-(1-|\nabla \rho|^2)\nabla\rho}{\rho^2}\\
    =&\ric_{\bar{g}}(\nabla\rho,\cdot),
\end{align*}
 where we used \eqref{Eq: scalar falt adapted metric - 3} in the second equality.

 On $\partial X$, we have 
 \[
 \nu=-\nabla \rho,
 \]
 which gives the expression of the second fundamental form of $\partial X$
 \begin{align}\label{Eq: second fundamental form and P}
     A=-\nabla^2\rho|_{T\partial X}=\frac{P}{2}h.
 \end{align}
 The boundary is therefore umbilic and the traced Gauss equation gives \eqref{Eq: scalar falt adapted metric - 6}.
\end{proof}

\begin{lem}\label{lem: local convexity of the boundary}
    Given a scalar flat adapted compactification $\bar{g}=\rho^2g_+$, for every two points in the interior of $B^4$, the shortest path connecting these two points is a geodesic in the interior of $B^4$.
\end{lem}
\begin{proof}
    This follows directly from the fact that the second fundamental form of $\partial X$ is positive.
\end{proof}

\subsection{Curvature bounds and a blow-up analysis}
We now apply a blow-up analysis to show the curvature bounds in Theorem \ref{Thm: compactness}. 
The main strategy follows from \cite{chang2025problem}*{\S 6} (see also \cite{changGe2018compactness}*{\S 4}).

\begin{proposition}\label{Prop: curvature bounds}
    Under the assumptions in Theorem \ref{Thm: compactness}, there exists a constant $C_1>0$ so that every adapted scalar flat compactification $\bar g_i$ of $g_+^i$ satisfies 
    \[\sup_{X}\left( |\Rm_{\bar g_i}| + |\nabla \Rm_{\bar g_i}|^{2/3} \right) \leq C_1.  \]
\end{proposition}

\begin{proof}
    Suppose the proposition fails. 
    Then there exists a family of conformally compact Einstein metrics $\{g_+^i\}$ on $B^4$ satisfying (i) and (ii) in Theorem \ref{Thm: compactness of chang-ge} so that the scalar flat adapted compactifications $\{\bar g_i=\rho_i^2g_+^i\}$ have
    \[ K_i^2 := \max_{\overline{B}^4}\left( |\Rm_{\bar g_i}|+ |\nabla\Rm_{\bar g_i}|^{2/3} \right)\to +\infty \]
    with the maximum attained at $p_i\in \overline{ B}^4$.  
    Set 
    \[ \hat g_i:= K_i^2 \bar g_i\qquad\text{and}\qquad \hat \rho_i:= K_i\rho_i. \]
    Then, $\hat g_i=\hat \rho_i^2\cdot g_+^i$ has 
    \begin{align}\label{Eq: normalized curvature bound}
        \max_{\overline{B}^4}\left( |\Rm_{\hat g_i}|+ |\nabla\Rm_{\hat g_i}|^{2/3} \right) = \left( |\Rm_{\hat g_i}|+ |\nabla\Rm_{\hat g_i}|^{2/3} \right)(p_i)  =1,
    \end{align}
    and Lemma \ref{Lem: computation for scalar flat adapted metric} remains valid for $(\hat g_i,\hat \rho_i)$ with 
    \[ \wha P_i=\frac{1-|d\hat\rho_i|^2_{\hat g_i}}{\hat\rho_i} = \frac{1-|d\rho_i|^2_{\bar g_i}}{K_i\rho_i} = K^{-1}_iP_i. \]

    We prove the following apriori bounds. 

\begin{lem}\label{lem:change-ge aprior bounds}
    For every $0\leq l\leq m-5$, there exists a constant $C
    _l$, independent of $i$, such that
    \begin{align}\label{eq: boudned l derivative of curvature}
        |\nabla^l\Rm_{\hat{g}_i}|\leq C_l, \quad \text{ on } X=\overline{B}^4,
    \end{align}
    and 
    \begin{align}\label{eq: bounded l derivative second fundamental form}
        |\nabla^l_{\partial X}A_i|\leq C_l, \quad \text{ on } \partial X=\bd\overline{B}^4,
    \end{align}
    where $A_i$ is the second fundamental form of $\partial X=\bd\overline{B}^4$, w.r.t. the metric $\hat{g}_i$.
\end{lem}
\begin{proof}
    The proof follows from applying the $\epsilon$ regularity result in \cite{chang2025problem}*{Corollary 3.3}.
    We need to check that for a fixed $i$, $(\overline{B}^4,\hat{g}_i)$ satisfies the conditions in \cite{chang2025problem}*{Corollary 3.3}, i.e. the assumptions in Theorem \ref{Thm: compactness of chang-ge} and the existence of a fixed scale $a\in(0,1]$  such that
    \begin{align}\label{eq: smallness of curvature}
        \|\Rm_{\hat{g}_i}\|_{L^2(B_a(p))}+\|H_i\|_{L^3(B_a(p)\cap \partial \overline{B}^4)}\leq \epsilon, \quad \forall p\in B^4, 
    \end{align}
    where $\epsilon\in(0,1)$ is a constant from \cite{chang2025problem}*{Corollary 3.3}.
    Note that assumptions (i) and (ii) in Theorem \ref{Thm: compactness of chang-ge} are already satisfied. In \cite{chang2025problem}*{Corollary 3.3}, it does not assume a uniform Weyl energy bound for a family of metrics, instead, it studies a fixed metric, whose Weyl energy is obviously bounded. Therefore, for a fixed index $i$, we only need to check that it satisfies the smallness hypothesis \eqref{eq: smallness of curvature}.

    Since we already have a uniform boundedness of $\Rm_{\hat{g}_i}$ and $\nabla \Rm_{\hat{g}_i}$, we only need the boundedness of $\vol_{\hat{g}_i}(B_a(p))$ and $\area_{\hat{g}_i}(B_a(p)\cap\partial\overline{B}^4)$.

    By Lemma \ref{lem: local convexity of the boundary} and the Bishop comparison theorem, since $\ric_{\hat{g}_i}\geq -3$, we have
    \[
    \vol_{\hat{g}_i}(B_a(p))\leq Ca^4.
    \]

    Let $\varphi$ be the standard cut-off function on $B_{2a}(p)$, such that $\varphi\equiv 1$ on $B_a(p)$, $\varphi\equiv 0$ outside $B_{2a}(p)$, and it has $|\nabla \varphi|\leq Ca^{-1}$. Using $\varphi$ as a test function, Equation \eqref{Eq: scalar falt adapted metric - 4} gives
    \[
    \int _{\partial X}\varphi \, da=-\int_X\langle \nabla \varphi,\nabla\hat{\rho}_i\rangle \, dV+2\int_X\varphi \wha{P}_i\, dV,
    \]
    where we used the fact that $\nabla\hat{\rho}_i=-\nu$ on $\partial \overline{B}^4$.

    On the boundary we have that $\wha{P}_i=\frac{2}{3}H_i$, which is uniformly bounded. By \eqref{Eq: scalar falt adapted metric - 5}, we know that for every point the $\nabla \wha{P}_i$ is uniformly bounded, hence we conclude that $\wha{P}_i$ is uniformly bounded by a linear function of its distance to the boundary. For $q\in \text{supp}(\varphi)$, we have $|\wha{P}_i|\leq C(1+a)$, where the constant $C$ does not depend on $a$ or the index $i$.

    With the volume estimates, it follows that 
    \[
    \area_{\hat{g}_i}(B_{a}(p)\cap \partial \overline{B}^4)\leq Ca^3.
    \]
    We can then find a small enough scale $a$ such that \eqref{eq: smallness of curvature} holds. Note that the scale of $a$ does not depend on the index of $i$, it comes from the universal curvature bound.

    Then \eqref{eq: boudned l derivative of curvature} follows directly from \cite{chang2025problem}*{Corollary 3.3}, and \eqref{eq: bounded l derivative second fundamental form} follows from \eqref{Eq: scalar falt adapted metric - 5} and the identity \eqref{Eq: second fundamental form and P}.
\end{proof}
    We now separate the proof into two cases.

    \vspace{1mm}
    {\bf Case I.} {\it Boundary blow-up.}

    Suppose $\dist_{\hat g_i}(p_i,\bd\overline{B}^4) \leq D<\infty$. 
    Then, take the boundary point $q_i\in\bd\overline{B}^4$ nearest to $p_i$. 
    Using \cite{chang2025problem}*{Lemma 6.1, 6.2} and Lemma \ref{lem:change-ge aprior bounds} in place of \cite{chang2025problem}*{Corollary 3.3}, the compactness argument of Chang-Ge in the proof of \cite{chang2025problem}*{Lemma 6.3} can be adapted to show that there exists a subsequence of $(\overline{B}^4,\hat g_i, q_i)$ converging in the pointed $C^{m-3,\beta}_{loc}$ ($\forall \beta\in (0,1)$) Cheeger–Gromov topology to a complete $C^{m-3,\beta}_{loc}$ pointed manifold $(N,g_\infty,q_\infty)$. 
    By the constructions, $(|\Rm_{g_\infty}|+ |\nabla\Rm_{g_\infty}|^{2/3})(p_\infty)=1$ at the limit point $p_\infty:=\lim p_i\in N$. 
    Additionally, the limit of $(S^3, K_i^2h_i, q_i)$ is the boundary $(\bd N, g_\infty|_{\bd N})$, which is a complete flat three-manifold with infinite injectivity radius (i.e., the Euclidean space $(\mb R^3,\delta)$) since $\{h_i\}$ is compact and $K_i\to+\infty$.

    Next, by passing to a subsequence, one can also obtain $\hat \rho_i\to\rho_\infty$ in $C^{m-3,\beta}_{loc}(N)$ and $\wha P_i\to P_\infty$ in $C^{m-4,\beta}_{loc}(N)$ so that \eqref{Eq: P quantity}-\eqref{Eq: scalar falt adapted metric - 5} are valid under the metric $g_\infty$. 
    Indeed, since $|\nabla\hat\rho_i|\leq 1$ and $\hat\rho_i|_{\bd \overline{B}^4}=0$, we have $0\leq \hat\rho_i \leq \dist_{\hat g_i}(\bd\overline{B}^4,\cdot)$. 
    Additionally, by \eqref{Eq: scalar falt adapted metric - 5}, \eqref{Eq: scalar falt adapted metric - 6}, and \eqref{Eq: normalized curvature bound}, we also have uniform bounds on $\wha P_i|_{\bd\overline{B}^4}=\frac{2}{3}H_{\hat g_i}$ and $|\nabla\wha P_i|\leq |\Ric_{\hat g_i}|$. 
    In particular, $\wha P_i$ is uniformly bounded on balls of fixed radius centered at $q_i$. 
    Together with \eqref{Eq: scalar falt adapted metric - 3}, we obtain a uniform bound of $\nabla^2\hat\rho_i$. 
    Using the $C^{m-3,\beta}_{loc}$-convergence of $(\overline{B}^4,\hat g_i,q_i)$, the higher derivatives of $\hat\rho_i$ and $\wha P_i$ can also be bounded by differentiating \eqref{Eq: scalar falt adapted metric - 3} and \eqref{Eq: scalar falt adapted metric - 5}. 
    Therefore, we obtain the convergences of $\hat \rho_i$ and $\wha P_i$ by Arzel{\`a}–Ascoli.

    The convergence implies that $\rho_\infty\geq 0$, $\rho_\infty|_{\bd N}=0$, $|\nabla\rho_\infty|\leq 1$, and $|\nabla\rho_\infty|= 1$ on $\bd N$ under the metric $g_\infty$. 
    If $\rho_\infty(x)=0$ for some interior point $x\in \interior(N)$, then $\nabla \rho_\infty(x)=0$, which contradicts $|d\rho_\infty|_{g_\infty}^2(x)=1-P_\infty(x)\rho_\infty(x)=1$. 
    This shows that $\rho_\infty^{-1}(0)=\bd N $, and thus $\rho_\infty$ is a defining function on $N$. 
    Moreover, one verifies that $g_+^\infty:=\rho_\infty^{-2}g_\infty$ is a complete metric in $\interior(N)$, which is also the $C^2_{loc}$-limit of $g_+^i=\hat\rho_i^{-2}\hat g_i$. 
    Therefore, $g_+^\infty$ is an Einstein metric with $\Ric_{g_+^\infty}=-3g_+^\infty$. 
    By the boundary regularity result (\cite{chrsciel2005BoundaryRegularity}), $g_+^\infty$ admits a {\em smooth} conformal compactification inducing the Euclidean boundary metric near each compact set of $\bd N$. 
    Hence, we see from the Liouville theorem (Theorem \ref{Thm: Liouville theorem for round filling}) that $g_+^\infty$ is hyperbolic. 
    Finally, by the proof of Chang-Ge in \cite{chang2025problem}*{Lemma 5.15, P. 55, Step 11}, we conclude that $g_\infty$ is the flat metric, which contradicts $(|\Rm_{g_\infty}|+ |\nabla\Rm_{g_\infty}|^{2/3})(p_\infty)=1$. 

    \vspace{1mm}
    {\bf Case II.} {\it Interior blow-up.}

    Suppose $\liminf_{i\to\infty}\dist_{\hat g_i}(p_i,\bd\overline{B}^4) = +\infty $. 
    Then, as in {\bf Case I}, one can adapt the compactness argument of Chang-Ge in the proof of \cite{chang2025problem}*{Lemma 6.3} to obtain a subsequence of $(\overline{B}^4, \hat g_i, p_i)$ that converges in the pointed $C^{m-3,\beta}_{loc}$ ($\forall \beta\in (0,1)$) Cheeger–Gromov topology to a complete $C^{m-3,\beta}_{loc}$ pointed manifold $(N,g_\infty,p_\infty)$ without boundary.
    Next, we make the following claim, which is similar to a statement in the proof of \cite{changGe2018compactness}*{Lemma 4.9}. 

    \begin{claim}\label{Claim: un-bounded value of defining function}
        Set $a_i:= \hat\rho_i(p_i)$. Then $\lim_{i\to\infty} a_i = +\infty$.
    \end{claim}
    \begin{proof}
        Suppose, after passing to a subsequence, $\sup_ia_i\leq a$ for some $a>0$. 
        Then, since $|\nabla \hat\rho_i|_{\hat g_i}=|\nabla \rho_i|_{\bar g_i}\leq 1$ in $B^4$, we have $\hat\rho_i|_{B_{R}(p_i)}\leq a+R$ for any $R$-ball in $(\overline{B}^4, \hat g_i)$. 
        
        Combining \eqref{Eq: adapted scalar-flat compactification} with \eqref{Eq: scalar falt adapted metric - 5} and \eqref{Eq: normalized curvature bound}, we conclude that $|d\wha P_i|_{\hat g_i} \leq |\Ric_{\hat g_i}|_{\hat g_i}\cdot |\nabla \hat \rho_i|_{\hat g_i}\leq C$ for some constant $C>0$. 
        Hence, if $\liminf_{i\to\infty}\hat P_i(p_i)=+\infty$, then for any arc-length parameterized geodesic $\gamma_i:[0,1]\to B_1(p_i)\subset (\overline{B}^4, \hat g_i)$ with $\gamma_i(0)=p_i$, we see from $|d\wha P_i|_{\hat g_i}\leq C$ that $\lim_{i\to\infty}\inf_{t\in[0,1]}\wha P_i(\gamma_i(t))=+\infty$ up to a subsequence. 
        Thus, by \eqref{Eq: scalar falt adapted metric - 3}, $(\hat\rho_i\circ\gamma_i)''=-\frac{\wha\rho_i}{2}\Ric_{\hat g_i}(\gamma_i',\gamma_i') - \frac{\wha P_i}{2}$ tends to $-\infty$ uniformly on $[0,1]$ as $i\to\infty$, which contradicts $|(\hat\rho_i\circ\gamma_i)'|\leq 1$ by integrating $(\hat\rho_i\circ\gamma_i)''$ over $[0,1]$. 
        Therefore, up to a subsequence, $\sup_i\wha P_i(p_i)\leq C$ for some $C>0$. 

        Now, as in {\bf Case I}, we can extract a subsequence so that $\hat \rho_i\to\rho_\infty \in C^{m-3,\beta}_{loc}(N)$ and $\wha P_i\to P_\infty\in C^{m-4,\beta}(N)$ with $\rho_\infty\geq 0$ on $N$ satisfying \eqref{Eq: P quantity}-\eqref{Eq: scalar falt adapted metric - 5}. 
        By the same arguments as in {\bf Case I}, we have $\rho_\infty>0$ in $N$. 
        Meanwhile, since $(N,g_\infty)$ is complete with $\Ric_{g_\infty}$ bounded from below, the Omori-Yau principle (cf. \cite{yau1975harmonic}) gives a sequence of points $\{x_k\}\subset N$ so that 
        \[ \rho_\infty(x_k)\to I:=\inf_{x\in N}\rho_\infty(x) \in [0,\infty),\qquad |\nabla\rho_\infty|_{g_\infty}(x_k)\to0,\qquad \text{and}\quad \Delta\rho_\infty(x_k)\geq -\frac{1}{k}.\]
        Therefore, for $k$ large enough, we have $|\nabla\rho_\infty|^2_{g_\infty}(x_k)<1/2$, and thus, by \eqref{Eq: P quantity},
        \[ \Delta\rho_\infty(x_k) = \frac{2(|\nabla\rho_\infty|^2_{g_\infty}(x_k)-1)}{\rho_\infty(x_k)} \leq -\frac{1}{\rho_\infty(x_k)}.  \]
        However, if $I>0$, then $\lim_{k\to\infty}\Delta\rho_\infty(x_k)\leq \lim_{k\to\infty}\frac{-1}{\rho_\infty(x_k)}=-\frac{1}{I}<0$; if $I=0$, then $\lim_{k\to\infty}\Delta\rho_\infty(x_k)\leq \lim_{k\to\infty}\frac{-1}{\rho_\infty(x_k)}=-\infty$. 
        This shows a contradiction. 
    \end{proof}

    Next, we can show that $g_\infty$ is Ricci flat (cf. \cite{changGe2018compactness}*{Lemma 4.9, Step 3}). 
    Specifically, consider the function $f_i:=\hat \rho_i/a_i$. 
    For fixed $R>0$, since $|\hat\rho_i-a_i|\leq R$ on $B_R(p_i)\subset (\overline{B}^4, \hat g_i)$, we have $f_i\to 1$ and $\nabla f_i\to 0$ uniformly in $B_R(p_i)$, which implies that $\nabla^2f_i$ tends to $0$ in the sense of distributions as $i\to\infty$. 
    By \eqref{Eq: scalar falt adapted metric - 3}, we also have 
    \[ \nabla^2 f_i = -\frac{1}{2}f_i\Ric_{\hat g_i}- \frac{\wha P_i}{2a_i}\hat g_i.  \]
    Therefore, since $0\leq \frac{\wha P_i}{a_i}\leq \frac{1}{a_i\hat\rho_i}\to 0$ as $i\to\infty$, we conclude that $\Ric_{g_\infty}=0$. 

    Finally, one can then apply the proof of Chang-Ge in \cite{changGe2018compactness}*{Proposition 4.15} to obtain a contradiction by topologically excluding such a limit $N$. 
    This proves the proposition. 
\end{proof}

\subsection{Proof of Theorem \ref{Thm: compactness}}

\begin{proof}[Proof of Theorem \ref{Thm: compactness}]
    Firstly, the curvature bound is given by Proposition \ref{Prop: curvature bounds}. 
    Combining with \eqref{Eq: scalar falt adapted metric - 6}, we see that the mean curvature $H_{\bar g_i}$ on $\bd\overline{B}^4$ is uniformly bounded. 
    Additionally, since the set $\{h_i\}$ of boundary metrics is compact, the boundary area $\Area_{h_i}(\bd\overline{B}^4)$ is also uniformly bounded. 
    Consequently, $\|H_{\bar g_i}\|_{L^1(\partial\overline B^4,h_i)}\leq C$ for some constant $C>0$. 
    Moreover, combining $Y(S^3,[h_i])\geq C_0>0$ with \cite{chang2025problem}*{Theorem 1.1 and (8.22)}, we conclude that $\vol_{\bar g_i}(\overline{B}^4)$ is uniformly bounded. 
    Finally, $\int_{\overline{B}^4}|W|^2_{\bar g_i} $ is uniformly bounded by the volume bound and the curvature bound in Proposition \ref{Prop: curvature bounds}. 
    This proves Theorem \ref{Thm: compactness}. 
\end{proof}

\section{Existence of equivariant Poincar{\' e}-Einstein fillings}

Given any $h\in Met_{\mb Z_p}(S^3)$ with $[h]\in \mc C_{\mb Z_p}$ satisfying $Y(S^3, [h])>0$, we show in this section that $\Pi^{-1}([h])\cap \mc E_{\mb Z_p}$ is non-empty. 

For any such metric $h\in Met_{\mb Z_p}(S^3)$, let $u_1>0$ with $| u_1|_{L^2}=1$ be the unique normalized first eigenfunction of the conformal Laplacian $L_h:=-\Delta_h+\frac{1}{8}R_h$. 
The uniqueness of $u_1$ and the $\langle\mk r_{p,q}\rangle$-invariance of $h$ imply that $u_1\circ \mk r_{p,q}=u_1$, and thus $u_1$ is a $\langle\mk r_{p,q}\rangle$-invariant function. 
Since $Y(S^3, [h])>0$, the first eigenvalue $\lambda_1$ of $L_h$ is positive. 
Thus, we have a $\langle\mk r_{p,q}\rangle$-invariant smooth metric $u_1^4h\in [h]$ with positive scalar curvature $R_{u_1^4h}>0$. 
In particular, we may assume without loss of generality that $(S^3,h)$ has positive scalar curvature $R_h>0$ by using $u_1^4h$ in place of $h$.

\subsection{Linearization of the Einstein operator with fixed gauge}

Let $\gamma\in Met_{\mb Z_p}(S^3)$ be a unit-density smooth boundary metric (cf. \eqref{Eq: unit-density representative}). 
Suppose $g_+$ is a smooth conformally compact Einstein metric on $B^4$ with conformal infinity $[\gamma]\in\mc C_{\mb Z_p}$. 
In this subsection, we investigate the linearization of the Einstein operator at $g_+$ with fixed gauge. 
In particular, we will show that $g_+$ is {\em weakly nondegenerate} in the sense of \cite{mazzeo2006maskit}*{Definition 4}. 
Combining this with the implicit function theorem, one can locally construct a finite dimensional $\mb Z_p$-invariant manifold in $\mc E$ containing $g_+$ (see the next subsection and \cite{mazzeo2006maskit}*{P. 395}).

The main results (Proposition \ref{Prop: gauge and chart} and \ref{Prop: linearization theory}) of this subsection are built on the local deformation theory of conformally compact Einstein metrics, the standard gauge fixing technique, and the analysis of the linearized Einstein operator on weighted H{\"o}lder spaces. 
We refer to \cite{anderson2008Einstein}\cite{anderson2010structure}\cite{biquard2000Metriques}\cite{lee2006fredholm}\cite{mazzeo2006maskit} for the development of these methods and related deformation results.
For completeness, we give detailed arguments explaining how
these results apply in the present setting and yield the precise
statements needed here.

\medskip
To begin with, we consider the deformations of the metric $g_+$ in the space of conformally compact metrics. 
Let $x$ be a fixed geodesic defining function on $B^4$ so that 
\begin{align}\label{Eq: CCE metric g_+ and fixed defining}
    g_+=x^{-2}(dx^2+ h_x)\qquad\text{with}\qquad h_x=\gamma+O(x^2). 
\end{align}
Then for any $m\geq 6$ and $\alpha,\delta\in (0,1)$, define the weighted H{\"o}lder space 
\[S^{m,\alpha}_\delta:=S^{m,\alpha}_\delta(B^4,g_+)\] 
as the space of symmetric bilinear forms $u$ on $(B^4,g_+)$ so that $\|x^{-\delta} u\|_{C^{m,\alpha}(B^4,g_+)}<\infty$ (cf. \cite{lee2006fredholm}*{\S 3}). 
Hence, for any $u$ in a small neighborhood of $0\in S^{m,\alpha}_\delta(B^4)$, 
\begin{align}\label{Eq: interior deformation}
    g_+\,+\, u
\end{align}
provides an interior deformation of the metric $g_+$. 
In particular, since $|x^2u|_{x^2g_+} = |u|_{g_+}=O(x^\delta)$, we have $\lim_{x\to 0}x^2(g_++u)=\gamma$ on $S^3$, which implies that such interior deformations of $g_+$ do not affect the boundary metric in the conformal compactification.

Next, for the `boundary deformation' of $g_+$, we use the
parameterization $F:\mc C\to\mc A$ from
Lemma \ref{Lem: parameterization of boundary conformal classes}. 
For each $A\in\mc A$, let $\gamma_A$ denote the unique
unit-density representative of $F^{-1}(A)\in\mc C$, given by
\[
    \gamma_A=h_0(e^A\cdot,\cdot),
    \qquad \text{with }\quad
    d\operatorname{Vol}_{\gamma_A}
    =d\operatorname{Vol}_{h_0}.
\]
Since $\gamma$ has the unit density, we can set
\[
    A_0:=F([\gamma])\qquad\text{so that } \qquad \gamma_{A_0}=\gamma.
\]
Then, for a neighborhood $\mc O_{A_0}\subset \mc A$ of $A_0$, $\{\gamma_A\}_{A\in \mc O_{A_0}}$ forms a smooth family of smooth metrics on $S^3$ with unit density. 
For a sufficiently small $\epsilon>0$, the collar $U_{2\epsilon}(S^3)=\{x<2\epsilon\}$ of $S^3$ in $B^4$ can be parameterized by $(0,2\epsilon)\times S^3$. 
Using the natural projection $\pi:(0,2\epsilon)\times S^3\to S^3$, we have a smooth family of symmetric $2$-tensors $\pi^*(\gamma_A-\gamma_{A_0})$ on the collar $U_{2\epsilon}(S^3)$. 
Fix a smooth cut-off function $\chi:[0,\infty)\to [0,1]$ so that $\chi |_{[0,\epsilon]}=1$ and $\chi|_{[2\epsilon,\infty)}=0$. 
Define then 
\begin{align}\label{Eq: boundary deformation}
    T_A:=g_+ + \chi(x)\cdot x^{-2}\cdot\pi^*(\gamma_A-\gamma_{A_0})= x^{-2}\left[ dx^2+ h_x+\chi(x)\cdot\pi^*(\gamma_A-\gamma_{A_0})  \right], 
\end{align}
which gives a boundary deformation of the metric $g_+=T_{A_0}$. 
In particular, for $\mc O_{A_0}\subset \mc A$ small enough, $T_A$ is a conformally compact metric on $B^4$ so that $x$ is a geodesic defining function of $T_A$ and the conformal compactification $x^2T_A$ has the boundary metric $\gamma_A$. 

Therefore, we have a family of conformally compact metrics defined by $T_A+u$, where $A\in \mc O_{A_0}$ and $u\in S^{m,\alpha}_\delta$ close to $0$.

\medskip
We next aim to find the Einstein metrics among these deformed metrics.  
For this purpose, define the Einstein operator 
\[ \Ein(g):= \Ric_g + 3g.\]
Then, every conformally compact Einstein deformation $T_A+u$ satisfies $\Ein(T_A+u)=0$. 
However, it is well-known that the linearization $D\Ein_{g_+}$ at the Einstein metric $g_+$ fails to be elliptic (cf. \cite{mazzeo2006maskit}*{\S 2.1.2}).  
Indeed, for any smooth vector field $X$ compactly supported in $B^4$ with the flow $\{\phi_t\}$, $\phi_t^*g_+$ remains an Einstein metric, i.e. $\Ein(\phi_t^*g_+)=\phi_t^*\Ein(g_+)=0$. 
Thus, the Lie derivative $\mc L_Xg_+$ satisfies $D\Ein_{g_+}(\mc L_Xg_+)=\left.\frac{d}{dt}\right|_{t=0}\Ein(\phi_t^*g_+)= 0$, which prevents $D\Ein_{g_+}$ from being elliptic by the arbitrariness of $X$ and the interior elliptic $H^2$ estimates. 
Nevertheless, for any metric $g'$ $C^{m,\alpha}$-close to $g_+$, define the Bianchi operator
\[ B_{g'}:S^{m,\alpha}_\delta\to C^{m-1,\alpha}_\delta(T^*B^4),\qquad   B_{g'}  (g)=\delta_{g'}g + \frac{1}{2}d \tr_{g'}g, \]
where $C^{l,\alpha}_\delta(T^*B^4)$ denotes the space of $1$-forms $\omega$ on $B^4$ with $|x^{-\delta}\omega|_{C^{l,\alpha}(B^4,g_+)}<\infty$, and $\delta_{g'}=-\div_{g'}$, i.e. $(\delta_{g'} u)_i=-\nabla_{g'}^j u_{ji}=-(g')^{jk}(\nabla^{g'}_ku)_{ji}$. 
Then by \cite{besse1987einstein}*{Theorem 1.174(d)}, 
\begin{align}\label{Eq: linearization of Ein}
    D\Ein_{g_+}(u) = \frac{1}{2} \Delta_L u -\delta^*\delta u - \frac{1}{2}\nabla d\tr u + 3u = Lu - \delta^* Bu, 
\end{align}
where we omitted the subscript $g_+$ for simplicity and used the following notations:
\begin{itemize}
    \item $\Delta_L$: the Lichnerowicz Laplacian,
    \item $\delta^*$: the symmetrized covariant derivative of $1$-forms, i.e. $\delta^*\omega=\frac{1}{2}(\nabla_i\omega_j + \nabla_j\omega_i)$, 
    \item $B:=B_{g_+}$ and $L=\frac{1}{2}\Delta_L+3 = \frac{1}{2}(\Delta_L+6) $. 
\end{itemize}
Therefore, if $Bu=0$, then $D\Ein_{g_+}(u) = Lu$ is elliptic. 
Moreover, it is more convenient to consider the following single operator. 

As in \cite{mazzeo2006maskit}*{(6)} (see also \cite{lee2006fredholm}*{(1.2)} and \cite{anderson2010structure}*{(4.18)}), for a sufficiently small $\mc O_{A_0}\subset \mc A$ and a small neighborhood $U$ of $0$ in $S^{m,\alpha}_\delta$, we can define a smooth map 
\begin{align}\label{Eq: Einstein operator with gauge}
    G : \mc O_{A_0}\times U &\longrightarrow S^{m-2,\alpha}_\delta , \nonumber
    \\
    (A,u)\quad &\longmapsto \Ein(T_{A}+u) + \delta_{T_{A}+u}^*B_{T_{A}}(T_{A}+u),
\end{align}
There are several benefits by considering the operator $G$. 
For instance, it follows from \cite{biquard2000Metriques}*{Lemma I.1.4} that any solutions of $G(A,u)=0$ give a conformally compact Einstein metric $T_A+u$ by shrinking $\mc O_{A_0}$ and $U$ if necessary. 
Indeed, we have the following result. 

\begin{proposition}\label{Prop: gauge and chart}
    Using the above notations, let $\mc Z(\mc O_{A_0}):= \wha\Pi^{-1}(\mc O_{A_0})$ be an open set of $\mc E$ equipped with the induced $C^\infty$ topology. 
    Then by shrinking the open sets $\mc O_{A_0}\subset\mc A$ and $ U\subset S^{m,\alpha}_\delta$ sufficiently small, the map 
    \[ \Psi: G^{-1}(0)\subset \mc O_{A_0}\times U \to \mc Z(\mc O_{A_0}), \qquad \Psi(A,u):= [T_A+u]_{PE},  \]
    is a well-defined homeomorphism onto an open neighborhood $\mc U_{z}\subset \mc Z(\mc O_{A_0})$ of $z:=[g_+]_{PE}$. 
\end{proposition}

\begin{remark}\label{Rem: regularity of deformations}
    If $T_A+u$ is a conformally compact Einstein metric with boundary metric $\gamma_A$ associated to the defining function $x$, then the boundary regularity results \cite{biquard2014analyse}*{Theorem 8.4} and \cite{chrsciel2005BoundaryRegularity}*{Theorem A} give a $C^{1,\alpha}$ collar diffeomorphism $\Phi$ on $\overline{ B}^4$ that is smooth in $B^4$ and $\Phi|_{\bd\overline{B}^4}=\Id_{S^3}$ so that $\Phi^*(T_A+u)$ is smoothly conformally compact. 
    Hence, we slightly abuse the notation and interpret $[T_A+u]_{PE}$ as $[\Phi^*(T_A+u)]_{PE}\in \mc E$ in the above proposition. 
\end{remark}

\begin{remark}\label{Rem: equivariant filling has equivariant deformation}
    If we further have $[g_+]_{PE}\in \mc E_{\mb Z_p}$, then by Proposition \ref{prop: action on the filling space}, there exists $\phi\in \Diff(\overline{B}^4)$ with $\phi|_{S^3}=\mk r_{p,q}$ so that $\phi^*g_+=g_+$ and $\phi$ generates a smooth $\mb Z_p$-action on $\overline{B}^4$. 
    Consequently, the geodesic defining function $x$ becomes $\langle\phi\rangle$-invariant and the boundary projection $\pi$ on the collar satisfies $\pi\circ\phi = \mk r_{p,q}\circ \pi$. 
    Additionally, since $F$ is $\langle \mk r_{p,q}\rangle$-equivariant, we conclude that $\mk r_{p,q}^* \gamma_A=\gamma_{\mk r_{p,q}\cdot A}$ and $\phi^* T_A=T_{\mk r_{p,q}\cdot A}$. 
    Therefore, the map 
    \[ \phi: (A,u)\longmapsto (\mk r_{p,q}\cdot A, \phi^*u) \]
    generates a $\mb Z_p$ action on $\mc A\times S^{m,\alpha}_\delta$ so that for any $(A,u)\in \mc O_{A_0}\times U$,
    \[ T_{\mk r_{p,q}\cdot A}+ \phi^*u= \phi^*(T_A+u)\qquad\text{and}\qquad G(\phi\cdot(A,u))=\phi^*G(A,u),  \]
    which implies that $G^{-1}(0)$ is $\mb Z_p$-invariant. 
    Moreover, recall that the map $\mk r_{p,q}$ on $S^3$ is radially extended to a diffeomorphism on $\overline{B}^4$ (still denoted by $\mk r_{p,q}$ in Proposition \ref{prop: action on the filling space}), which generates the $\mb Z_p$-action on $\mc E$. 
    Using the extended $\mk r_{p,q}$, $\eta:=\mk r_{p,q}^{-1}\circ\phi \in \Diff_{\bd}(\overline{B}^4)$ and 
    \begin{align*}
        \Psi(\phi\cdot (A,u)) &=[T_{\mk r_{p,q}\cdot A} + \phi^*u]_{PE} =[\phi^*(T_A+u)]_{PE} = [\eta^*(\mk r_{p,q}^*(T_A+u))]_{PE}
        \\
        &=[\mk r_{p,q}^*(T_A+u)]_{PE}=\mk r_{p,q}\cdot [T_A+u]_{PE} = \mk r_{p,q}\cdot \Psi(A,u),
    \end{align*}
    where we used $\mc E=E_{PE}/\Diff_{\bd}(\overline{B}^4)$ in the second line. This shows that $\Psi$ is $\mb Z_p$-equivariant when $[g_+]_{PE}\in \mc E_{\mb Z_p}$.
\end{remark}

\begin{proof}
    For any $(A,u)$ with $G(A,u)=0$, we see from $\Ric_{g_+}=-3g_+$ that $\Ric_{T_A+u}<0$ by shrinking $\mc O_{A_0}$ and $ U$. 
    Additionally, since $B_{T_A}(T_A)=0$, we have $\omega:=B_{T_A}(T_A+u)=B_{T_A}(u)\in C^{m-1,\alpha}_\delta(T^*B^4) $ and $|\omega|_{T_A+u}=O(x^\delta)\to 0$ as $x\to 0$. 
    Hence, by \cite{biquard2000Metriques}*{Lemma I.1.4} (see also \cite{mazzeo2006maskit}*{Proposition 1}), $T_A+u$ is an Einstein metric which is in the Bianchi gauge relative to $T_A$, i.e. $B_{T_A}(u)=B_{T_A}(T_A+u)=0$. 
    In particular, combining with the boundary regularity \cite{chrsciel2005BoundaryRegularity} and Remark \ref{Rem: regularity of deformations}, we have $[T_A+u]_{PE}\in \mc Z (\mc O_{A_0})$ and $\Psi$ is well defined. 

    Let $Z\in C_{\delta}^{m+1,\alpha}$ be any vector field, and define $\varphi_Z(p):=\exp_p^{g_+}(Z(p))$. Then it is a diffeomorphism of $B^4$ provided $\|Z\|_{C^{m+1,\alpha}_{\delta}}$ is sufficiently small. Note that
    \begin{align*}
    B_{T_A}((\varphi_Z^{-1})^*(T_A+u))=(\varphi^{-1}_Z)^*B_{\varphi_Z^*T_A}(T_A+u)=(\varphi_Z^{-1})^*B_{\varphi_Z^*T_A}(T_A+u-\varphi_Z^*T_A),
    \end{align*}
    where we used the fact that $B_{T_A}(T_A)=0$.
    It is straightforward from the definition that the following function
    \[
    F(A,u,Z)=B_{\varphi_Z^*T_A}(T_A+u-\varphi_Z^*T_A)
    \]
    is a smooth map from $\mathcal O_{A_0}\times S_{\delta}^{m,\alpha}\times C_{\delta}^{m+1,\alpha}\to C_{\delta}^{m-1,\alpha}(T^*B^4)$.
    Direct computations give
    \[
    D_ZF(A_0,0,0)Y=-(\nabla^*\nabla+3)Y^\flat,
    \]
    which is Fredholm of index 0, hence it is an isomorphism. The implicit function theorem then gives that $Z(A,u)$, defined as $F(A,u,Z)=0$, is a continuous map, provided $\mathcal O_{A_0}$ and $U$ are sufficiently small. If $B_{T_A}(u)=0$, we have $Z(A,u)=0$ by uniqueness. By shrinking $\mathcal O_{A_0}$ and $U$, the injectivity of $\Psi$ follows easily from the uniqueness of $Z$.

    We next prove $\Psi$ is continuous. Suppose $(A_i,u_i)\to (A,u)$ in $G^{-1}(0)$, by Theorem \ref{Thm: compactness}, $\Psi(A_i,u_i)$ converges (after passing to a subsequence) in $\mathcal E$, denoted as $z_{\infty}$. Take representatives $g_i\to g_{\infty}$ of $\Psi(A_i,u_i)$ and $z_{\infty}$, and isometries $\theta_i|_{S^3}=\text{Id}$, $\theta_i^*(T_{A_i}+u_i)=g_i$. $\theta_i$ are equicontinuous on $\overline{B^4}$. By the Arzela-Ascoli theorem, a subsequence converges to an isometry $\theta$ with $\theta|_{S^3}=\text{Id}$, $\theta^{*}(T_A+u)=g_{\infty}$, hence $z_{\infty}=[\theta^*(T_A+u)]_{PE}=\Psi(A,u)$.

    Let $D=G^{-1}(0)\cap (\mathcal O_{A_0}\times U)$, if $\Psi(D)$ does not contain an open neighborhood of $[g_+]_{PE}$, then there exists $g_i\to g_+$ in $ E_{PE}$ with $[g_{i}]_{PE}\notin \Psi(D)$. Let $A_i=\wha \Pi([g_i]_{PE})$, it is straightforward to show that $\tilde{u}_i:=g_i-T_{A_i}\to 0\in S_{\delta}^{m,\alpha}$. Then we have $Z_i:=Z(A_i,\tilde{u}_i)\to 0$, and $u_i=(\varphi_{Z_i}^{-1})^*g_i-T_{A_i}\to 0\in S_{\delta}^{m,\alpha}$, $G(A_i,u_i)=0$. Hence we have $[g_i]_{PE}=\Psi(A_i,u_i)\in\Psi(D)$ which is a contradiction.

    It then follows easily that $\Psi$ is indeed a well-defined homeomorphism onto the neighborhood $\mathcal U_z$.
\end{proof}

By the above proposition, a small neighborhood $\mc U_z$ of $z:=[g_+]_{PE}$ in $\mc E$ can be parameterized by $G^{-1}(0)\subset \mc O_{A_0}\times  U$. 
In the following proposition, we further show that $g_+$ is weakly nondegenerate in the sense of \cite{mazzeo2006maskit}*{Definition 4}. 

\begin{proposition}\label{Prop: linearization theory}
    The tangent map $DG_{(A_0,0)}: \mc A\times S^{m,\alpha}_\delta  \to S^{m-2,\alpha}_\delta $ is given by 
    \[ DG_{(A_0,0)}(v,w)= D_AG_{(A_0,0)}(v)+ D_uG_{(A_0,0)}(w) \]
    so that the following statements hold. 
    \begin{itemize}
        \item[(i)] $D_uG_{(A_0,0)}=L$ is a Fredholm linear map with Fredholm index $0$. In particular, the kernel $\mc N:=\ker L$ is of finite dimension, and the range $\im L$ is $\mc N^\perp:=\{f: \int_{B^4} \langle f,\kappa\rangle_{g_+} dV_{g_+} =0 \mbox{ for all }\kappa\in\mc N \}$.
        \item[(ii)] Every $\kappa\in\mc N$ is smooth, $g_+$-trace-free, and divergence-free (i.e. $\tr_{g_+}\kappa=0$ and $\delta_{g_+}\kappa=0$). Additionally, using the fixed geodesic defining function $x$, $\kappa\in \mc N$ has the expansion 
        \[ \kappa= x^{-2}\left[x^3\kappa_{(3)}+ O(x^4\log(x))\right] \qquad\mbox{near $S^3=\{x=0\}$ in $B^4$}\] 
        so that $\kappa_{(3)}$ is a smooth, $\gamma$-trace-free, symmetric $2$-tensor on $S^3$. 
        \item[(iii)] For any $\kappa\in\mc N$ and $v\in \mc A$, 
            \begin{align}\label{Eq: boundary realization}
                \int_{B^4} \langle D_AG_{(A_0,0)}(v),\kappa\rangle_{g_+} \, dV_{g_+}= -\frac{3}{2}\int_{S^3}\langle J_{A_0}(v) , \kappa_{(3)} \rangle_{\gamma} \,dA_\gamma, 
            \end{align}
            where $J_{A_0}(v):=\frac{d} {dt}\big|_{t=0}\gamma_{A_0+tv}$ is a smooth  $\gamma_{A_0}$-trace-free symmetric $2$-tensor on $S^3$.
        \item[(iv)] If $\kappa\in\mc N$ has $\kappa_{(3)}=0$, then $\kappa=0$. 
    \end{itemize}
\end{proposition}
\begin{proof}
    {\bf Step 1.}
    Firstly, the statements in (i) essentially follow from \cite{lee2006fredholm}*{Theorem C and 6.6}, see also \cite{mazzeo2006maskit}*{Proposition 4 and 5, and the statement after (16)}. 
    For completeness, we include more details. 
    Since $\Ein(g_+)=0$ and $B_{g_+}(g_+)=0$, we see from \eqref{Eq: linearization of Ein} and \eqref{Eq: Einstein operator with gauge} that 
    \[D_uG_{(A_0,0)}=D\Ein_{g_+} + \delta^*B = L = \frac{1}{2}(\Delta_L+6).\] 
    By \cite{lee2006fredholm}*{Proposition D} with $n=3$ and $c=6$, the operator $L$ satisfies the hypotheses of \cite{lee2006fredholm}*{Theorem C(c)} with $n=3$ and $R=\frac{3}{2}$. 
    Hence, \cite{lee2006fredholm}*{Theorem C(c)} implies that for any $\delta'\in (0,3)$, $L: S^{m,\alpha}_{\delta'}\to S^{m-2,\alpha}_{\delta'}$ is a Fredholm map of Fredholm index $0$. 
    Since we take $\delta\in (0,1)$, the first statement in (i) is valid. 
    In particular, $\dim(\mc N)<\infty$ and 
    \begin{align}\label{Eq: decay of kappa in kernel}
        \mc N\subset S^{m,\alpha}_{\delta'},\qquad \forall\delta'\in (0,3).
    \end{align} 
    Thus, for any $\kappa\in\mc N$ and $f\in S^{m-2,\alpha}_\delta$ with $\delta\in (0,1)$, we can take $\delta'\in (3-\delta, 3)$ so that on the collar $\{0<x<\epsilon\}\cong (0,\epsilon)\times S^3$ of $S^3$ in $B^4$, $|f|_{g_+}=O(x^\delta)$, $|\kappa|_{g_+}=O(x^{\delta'})$, $dV_{g_+}=x^{-4}\,dx\,d\sigma_{h_x}$, and
    \begin{align}\label{Eq: decay estimates for kernel}
        |\langle f, \kappa\rangle| \,dV_{g_+} \leq C x^{\delta+\delta'-4} \,dx\,d\sigma_{h_x},
    \end{align}
    where $C>0$ is a constant and $h_x$ is given in \eqref{Eq: CCE metric g_+ and fixed defining}. 
    Since $\delta+\delta'-4>-1$, $\langle f, \kappa\rangle$ is $dV_{g_+}$-integrable near the boundary $\{x=0\}$. 
    Now, for any $f\in \im L$ with $f=Lw$, integrating by parts in the compact domain $\Omega_{\tilde\epsilon}:=\{x\geq \tilde\epsilon\}$ shows that 
    \begin{align*}
        \int_{\Omega_{\tilde\epsilon}}\langle L w, \kappa\rangle_{g_+} - \langle w, L\kappa\rangle_{g_+} dV_{g_+} &= \frac{1}{2}\int_{\bd\Omega_{\tilde\epsilon}} \langle w, \nabla_{\nu_{\tilde\epsilon}} \kappa\rangle_{g_+} - \langle \nabla_{\nu_{\tilde\epsilon}} w,  \kappa\rangle_{g_+} d\sigma_{g_+|_{\{x=\tilde\epsilon\}}}
        \\
        &\leq C \int_{S^3} \tilde\epsilon^{\delta+\delta'-3} d\sigma_{h_{\tilde\epsilon}} \to 0 \qquad\mbox{as $\tilde\epsilon\to 0$},
    \end{align*}
    where $\nu_{\tilde\epsilon}=-x\bd_x$ is the outward pointing unit normal, and we used $|w|+|\nabla w|=O(x^\delta)$, $|\kappa|+|\nabla \kappa|=O(x^{\delta'})$, and $d\sigma_{g_+|_{\{x=\tilde\epsilon\}}}=x^{-3}\, d \sigma_{h_{\tilde\epsilon}}$ in the inequality. 
    Therefore, 
    \begin{align*}
        \int_{B^4}\langle f, \kappa\rangle_{g_+}dV_{g_+}= \int_{B^4}\langle Lw, \kappa\rangle_{g_+}dV_{g_+}=\int_{B^4}\langle w, L\kappa\rangle_{g_+}dV_{g_+}=0,
    \end{align*}
    which implies $\im L \subset \mc N^\perp$. 
    By the zero Fredholm index of $L$, $\codim(\im L )=\codim(\mc N^\perp)=\dim(\mc N)$. 
    This shows $\im L  = \mc N^\perp$. 

    \medskip
    {\bf Step 2.}
    Next, for any $\kappa\in\mc N$, since $\tr_{g_+}(\Delta_L )=\nabla^*\nabla\tr_{g_+}$, taking the trace of $L\kappa=0$ gives $(\nabla^*\nabla+6)\tr_{g_+}\kappa=0$. 
    By the decay estimates as in \eqref{Eq: decay of kappa in kernel} and \eqref{Eq: decay estimates for kernel}, one can integrate $\langle (\nabla^*\nabla+6)\tr_{g_+}\kappa, \tr_{g_+}\kappa\rangle_{g_+} $ by parts on $\Omega_{\tilde\epsilon}=\{x\geq \tilde\epsilon\}$ and take $\tilde\epsilon\to 0$ to show that 
    \[ \int_{B^4} |d(\tr_{g_+}\kappa)|^2_{g_+} + 6 |\tr_{g_+}\kappa|^2 \, dV_{g_+} = 0, \]
    which implies that $\tr_{g_+}\kappa=0$. 
    Similarly, since $B\circ D\Ein_{g_+}=0$ and $2B\delta^*\omega=\nabla^*\nabla \omega-\Ric_{g_+}(\omega) = (\nabla^*\nabla+3)\omega$ for any $1$-form $\omega$, applying $B$ to $L\kappa=0$ gives $(\nabla^*\nabla+3)B\kappa=0$. 
    Combining the decay estimates \eqref{Eq: decay of kappa in kernel} and \eqref{Eq: decay estimates for kernel} with the integration by parts on $\Omega_{\tilde\epsilon}$ with $\tilde\epsilon\to 0$, we also conclude that 
    \[ \int_{B^4} |\nabla B\kappa|^2_{g_+} + 3|B\kappa|^2_{g_+} dV_{g_+} = 0, \]
    which implies $0=B\kappa=\delta_{g_+}\kappa+\frac{1}{2}d\tr_{g_+}\kappa=\delta_{g_+}\kappa$. 

    Additionally, the expansion in (ii) can be obtained from \cite{mazzeo2006maskit}*{(13)}. 
    Specifically, recall that $x^2g_+$ is smooth in $\overline{B}^4$ (by the boundary regularity \cite{chrsciel2005BoundaryRegularity}). 
    Using the fixed geodesic defining function $x$ on $B^4$, we have $dV_{g_+}= x^{-4}\,dx\,d\sigma_{h_x}$ by \eqref{Eq: CCE metric g_+ and fixed defining}. 
    Consider the density $dV_b:= x^{-1}\, dx\, d\sigma_{h_x}= x^{3}\, dV_{g_+}$ on $B^4$ (cf. \cite{mazzeo2006maskit}*{\S 2.2.2}). 
    Since the decay estimates in \eqref{Eq: decay estimates for kernel} provide the finite $L^2(B^4,g_+)$ norm of any $\kappa\in \mc N$, we have 
    \[\|x^{-3/2}\kappa\|_{L^2_b}^2 := \int_{B^4} x^{-3}|\kappa|_{g_+}^2 dV_b =\int_{B^4}  |\kappa|_{g_+}^2 dV_{g_+}=\|\kappa\|_{L^2(B^4,g_+)}^2<\infty,\]
    which implies that $\kappa\in x^{\frac{3}{2}}L^2_b$ satisfies the hypotheses of \cite{mazzeo2006maskit}*{Proposition 5} with $\mu:=\frac{3}{2}$ (see \cite{mazzeo2006maskit}*{P. 391} for the notations). 
    Therefore, \cite{mazzeo2006maskit}*{(13)} gives an expansion of $\kappa$: 
    \[\kappa \sim \sum \kappa_{j,l,r}\cdot x^{\zeta_j^++l}(\log x)^r\qquad\mbox{near the boundary $\{x=0\}$},\]
    where 
    \begin{itemize}
        \item $\kappa_{j,l,r}$ locally has the form of $\sum_{a,b=0}^3 f^{j,l,r}_{ab} \frac{dw_a}{x}\otimes\frac{dw_b}{x}$ for some smooth function $f^{j,l,r}_{ab}$ on $\{x=0\}$ under the local coordinates $(w_0=x,w_1,w_2,w_3)$ of collars,
        \item $\{\zeta_j^+\}_{j=1}^3$ are the indicial roots of $L$ for which $x^{\zeta_j^+}\in x^{\frac{3}{2}}L^2_b$. 
    \end{itemize}
    By \cite{mazzeo2006maskit}*{Proposition 3} with $n=3$, the smallest indicial root is $\zeta_3^+=3$. 
    One can also see from \eqref{Eq: decay of kappa in kernel} that the lowest term is $\kappa_{(3)}x^3=\kappa_{3,0,0}x^3$ in the above expansion. 
    Additionally, \cite{mazzeo2006maskit}*{Proposition 5} also shows that $\kappa\in O(x^3)$. 
    In particular, there is no $x^3(\log x)^r$ term for any $r\geq 1$ in the above expansion, and thus $\kappa\sim \kappa_{(3)}x^3+ O(x^4\log(x))$ near the boundary $\{x=0\}$. 
    Applying the computation and the explanation of \cite{mazzeo2006maskit}*{(11)} with $h:=\kappa$ and $(s,u_i):=(x,w_i)$, we conclude that $\kappa_{(3)}$ only has the tangential part, i.e. $f^{3,0,0}_{a,0}=f^{3,0,0}_{0,b}=0$ for $a,b\in\{0,\dots,3\}$, which implies that $\kappa_{(3)}$ is a smooth $\gamma$-trace-free symmetric $2$-tensor on $S^3$. 
    Finally, since the standard basis $\{dw_i\otimes dw_j\}$ of symmetric $2$-tensors is rescaled by $x^{-2}$ in \cite{mazzeo2006maskit}*{P. 388}, we obtain the expansion in (ii). 

    \medskip
    {\bf Step 3.}
    At $u=0$, we have $B_{T_A}(T_A)=0$ and $G(A,0)=\Ein(T_A)$. 
    Hence, given $v\in\mc A$ with $J_{A_0}(v):=\frac{d }{dt}\big|_{t=0} \gamma_{A_0+ tv}$, we have $\dot g:= \frac{d}{d t}\big|_{t=0} T_{A_0+tv}=\chi(x)x^{-2} \pi^*(J_{A_0}(v))$, and by \eqref{Eq: linearization of Ein},
    \[ D_AG_{(A_0,0)}(v) = D\Ein_{g_+}(\dot g) = L\dot g - \delta_{g_+}^* B\dot g .\]
    For any $\tilde\epsilon>0$, integrating by parts on $\Omega_{\tilde\epsilon}:=\{x\geq \tilde\epsilon\}$ shows that
    \begin{align*}
        I_{\tilde\epsilon} 
        &:= \int_{\Omega_{\tilde\epsilon}} \langle D_AG_{(A_0,0)}(v),\kappa\rangle_{g_+}\,dV_{g_+}
        =\int_{\Omega_{\tilde\epsilon} }\langle L\dot g,\kappa\rangle_{g_+}\,dV_{g_+} - \int_{\Omega_{\tilde\epsilon}} \langle \delta^*_{g_+}B\dot g, \kappa\rangle_{g_+}\,dV_{g_+}
        \\
        &= \int_{\Omega_{\tilde\epsilon}} \langle \dot g, L\kappa\rangle_{g_+} - \langle B\dot g, \delta_{g_+}\kappa\rangle_{g_+} \,dV_{g_+}      
        \\
        &\qquad +  \int_{\bd\Omega_{\tilde\epsilon}} \frac{1}{2}\langle \dot g,\nabla_\nu\kappa\rangle_{g_+} - \frac{1}{2}\langle \nabla_\nu \dot g, \kappa\rangle_{g_+} - \langle B\dot g, \kappa(\nu)\rangle_{g_+} \,d\sigma_{g_+|_{\bd\Omega_{\tilde\epsilon}}},
    \end{align*} 
    where $\nu=-x\bd_x$ is the outward pointing unit normal, $L=\frac{1}{2}\Delta_L+3$, and $\delta_{g_+}=-\div_{g_+}$. 
    Since $\kappa\in\mc N=\ker L $ is divergence free (see (ii)), we conclude that 
    \begin{align}\label{Eq: integration by parts for D_sG}
        \int_{\Omega_{\tilde\epsilon}} \langle D_AG_{(A_0,0)}(v),\kappa\rangle\,dV_{g_+} = \int_{\bd\Omega_{\tilde\epsilon}} \frac{1}{2}\langle \dot g,\nabla_\nu\kappa\rangle - \frac{1}{2}\langle \nabla_\nu \dot g, \kappa\rangle - \langle B\dot g, \kappa(\nu)\rangle \,d\sigma_{g_+|_{\bd\Omega_{\tilde\epsilon}}}.
    \end{align}
    (See also \cite{anderson2010structure}*{(4.30)} and \cite{mazzeo2006maskit}*{(18)}.)

    For the first boundary integration, take a local coordinate system $(w_0=x,w_1,w_2,w_3)$ on a collar, and write 
    \[ \dot g_{ij} = x^{-2} \dot \gamma_{ij}, \qquad \kappa_{ij}=x(\kappa_{(3)})_{ij}+O(x^2\log x). \]
    Then, under the metric $g_+=x^{-2}(dx^2+h_x)$ with $h_x=\gamma+ O(x^2)$, we have $\Gamma_{0,i}^k= - \frac{1}{x}\delta^k_i + O(x)$, and thus
    \[(\nabla_{\nu}\kappa)_{ij}=-x(\bd_x\kappa_{ij}-\Gamma^k_{0i}\kappa_{kj}-\Gamma^k_{0j}\kappa_{ik})=-3x(\kappa_{(3)})_{ij}+ O(x^2\log x).\]
    Hence, $\langle \dot g, \nabla_{\nu}\kappa\rangle_{g_+}=-3 x^3\langle J_{A_0}(v), \kappa_{(3)}\rangle_{\gamma} + O(x^4\log x)$.
    Combining with $d\sigma_{g_+|_{\bd\Omega_{\tilde\epsilon}}}=x^{-3}\,d\sigma_{h_x}$, we conclude that 
    \begin{align}\label{Eq: boundary integration 1}
        \int_{\bd\Omega_{\tilde\epsilon}} \frac{1}{2}\langle \dot g,\nabla_\nu\kappa\rangle_{g_+}\, d\sigma_{g_+|_{\bd\Omega_{\tilde\epsilon}}} = -\frac{3}{2}\int_{S^3} \langle J_{A_0}(v), \kappa_{(3)}\rangle_{\gamma}\, d\sigma_\gamma + O(\tilde\epsilon \log \tilde\epsilon)   .
    \end{align}
    By a similar estimate, we have $|\nabla_\nu\dot g|_{g_+}=O(x^2)$. 
    Combining with $|\kappa|_{g_+}=O(x^3)$, we conclude that 
    \begin{align}\label{Eq: boundary integration 2}
        \left| \int_{\bd \Omega_{\tilde\epsilon}} \langle \nabla_\nu \dot g, \kappa\rangle_{g_+}\,d\sigma_{g_+|_{\bd\Omega_{\tilde\epsilon}}} \right| = O(\tilde \epsilon^2). 
    \end{align}
    Moreover, since $\kappa_{(3)}$ is a $2$-tensor on $S^3$ that only has the tangential terms, we see that $\kappa_{(3)}(\nu,\cdot )=0$, and thus $|\kappa(\nu,\cdot)|_{g_+}=O(x^4\log x)$ by the expansion in (ii). 
    Using $\dot g=\chi \cdot (\pi^*J_{A_0}(v))_{ij}\frac{dw_i}{x}\otimes\frac{dw_j}{x}$, one verifies that $|B\dot g|_{g_+}=O(1)$. 
    This shows 
    \begin{align}\label{Eq: boundary integration 3}
        \left| \int_{\bd\Omega_{\tilde\epsilon}}  \langle B\dot g, \kappa(\nu)\rangle_{g_+} \,d\sigma_{g_+|_{\bd\Omega_{\tilde\epsilon}}} \right| = O(\tilde\epsilon \log \tilde\epsilon) .
    \end{align}
    Combining \eqref{Eq: integration by parts for D_sG} with \eqref{Eq: boundary integration 1}-\eqref{Eq: boundary integration 3}, we obtain \eqref{Eq: boundary realization} by letting $\tilde\epsilon\to 0$. 

    \medskip
    {\bf Step 4.}
    Suppose $\kappa_{(3)}=0$. 
    Then, $|\kappa|_{g_+}=|x^2\kappa|_{x^2g_+}=O(x^4\log(x))=o(x^3)$ near the boundary $\{x=0\}$. 
    If $\kappa$ is in the geodesic gauge, i.e. $\partial_{x}\lrcorner \kappa:= \kappa(\bd_x,\cdot) =0$, then by \cite{biquard2008continuation}*{Theorem 2}, we directly obtain $\kappa=0$. 
    
    If $\partial_{x}\lrcorner \kappa=\kappa(\bd_x,\cdot)\neq 0$, then we can first put $\kappa$ in the geodesic gauge by an infinitesimal diffeomorphism. 
    Namely, there is a smooth vector field $X$ on $B^4$ so that $\tilde\kappa:= \kappa- \mc L_Xg_+$ satisfies $\partial_{x}\lrcorner\tilde\kappa=\tilde\kappa(\bd_x,\cdot)=0$ near $\{x=0\}$. 
    Indeed, write $X=a(x,y)\bd_x + b^i(x,y)\bd_{y_i}$ on a collar $\{0<x<\epsilon\}$ with coordinates $(x,y_1,y_2,y_3)$. 
    We can solve the equation 
    \[ 
    \partial_{x}\lrcorner\tilde\kappa=0\Leftrightarrow \left\{
    \begin{aligned}
        0 &= \tilde\kappa_{xx}= \kappa_{xx}-(\mc L_Xg_+)_{xx} = \kappa_{xx}- 2x^{-2}\left(\bd_xa-\frac{a}{x}  \right)=  \kappa_{xx}- 2x^{-1}\bd_x\left(\frac{a}{x}\right)
        \\
        0 &= \tilde\kappa_{xi}=\kappa_{xi}-(\mc L_Xg_+)_{xi} = \kappa_{xi} - x^{-2}\left( \bd_ia + (h_x)_{ij}\bd_x b^j \right)
    \end{aligned}
    \right.
    \]
    with boundary condition $\lim_{x\to0}\frac{a(x,y)}{x}= b^i(0,y)=0$ to obtain $a(x,y)=\frac{x}{2}\int_0^x s\cdot\kappa_{xx}(s,y) \,ds$ and $b^i(x,y)=\int_0^x(h_s)^{ij} \left(s^2\kappa_{xj}(s,y)-\partial_j a(s,y)\right)\,ds$. 
    This defines $X$ on the collar $T_\epsilon:=\{0<x<\epsilon\}$ so that $\tilde\kappa=\kappa - \mc L_Xg_+$ is in the geodesic gauge. 

    By the construction of $X$ on $T_\epsilon$, we have 
    \begin{align}
        |X|_{g_+} + |\nabla X|_{g_+} = O(x^4\log(x)) \qquad\text{on }T_\epsilon. 
    \end{align}
    Consequently, $\tilde\kappa=\kappa- \mc L_Xg_+$ satisfies $|\tilde\kappa|_{g_+} = O(x^4\log(x)) = o(x^3)$ on the collar $T_\epsilon$. 
    We can now apply \cite{biquard2008continuation}*{Theorem 2, 3} to conclude that $0=\tilde\kappa = \kappa- \mc L_Xg_+$ on the collar $T_\epsilon$. 
    
    Additionally, by \cite{deturck1981regularity}*{Theorem 5.2}, the Einstein metric $g_+$ is analytic in the interior harmonic and geodesic normal coordinates, which implies the interior analyticity of $\kappa$ with $L_{g_+}\kappa=0$ by the elliptic regularity. 
    Since $\pi_1(B^4,T_\epsilon)=0$, it then follows from  \cite{An2024static}*{Theorem 9} that the $\kappa$-Killing field $X$ (i.e. a vector field with $\mc L_Xg_+=\kappa$) on the collar $T_\epsilon$ can be extended to a unique $\kappa$-Killing field on $B^4$. 
    This gives a global vector field $X$ on $B^4$ with $\kappa=\mc L_Xg_+$. 

    Finally, note that $\kappa=\mc L_Xg_+=2\delta^*_{g_+}X^\flat$, where $X^\flat$ is the $g_+$-metric dual $1$-form of $X$. 
    Integrating by parts on $\Omega_{\tilde\epsilon}=\{x>\tilde\epsilon\}$ shows that 
    \begin{align*}
        \int_{\Omega_{\tilde\epsilon}} |\kappa|_{g_+}^2 \,dV_{g_+} &= 2 \int_{\Omega_{\tilde\epsilon}} \langle \kappa, \delta^*_{g_+}X^\flat\rangle \, dV_{g_+} = 2\int_{\Omega_{\tilde\epsilon}} \langle \delta_{g_+}\kappa, X^\flat\rangle \, dV_{g_+} + 2\int_{\bd \Omega_{\tilde\epsilon}} \kappa(\nu, X) d\sigma_{g_+|_{\bd\Omega_{\tilde\epsilon}}}
        \\
        &= 2\int_{\bd \Omega_{\tilde\epsilon}} \kappa(\nu, X) d\sigma_{g_+|_{\bd\Omega_{\tilde\epsilon}}},
    \end{align*}
    where $\nu$ is the outward pointing unit normal, and $\delta_{g_+}\kappa=0$ (cf. (ii)) is used in the second line. 
    Combining $|\kappa|_{g_+}=O(x^4\log(x))$ with $|X|_{g_+}=O(x^4\log(x))$ and $d\sigma_{g_+|_{\bd\Omega_{\tilde\epsilon}}}=x^{-3}\,d\sigma_{h_x}$, we have $|\kappa(\nu,X)|d\sigma_{g_+|_{\bd\Omega_{\tilde\epsilon}}} \leq C \tilde\epsilon^5|\log(\tilde\epsilon)|^2 d\sigma_{h_x} $ for some constant $C>0$. 
    By taking $\tilde\epsilon\to 0$, we conclude that  $\int_{B^4} |\kappa|_{g_+}^2 \,dV_{g_+} = 0$ and $\kappa=0$. 
    This proves (iv).
\end{proof}

\subsection{A finite dimensional equivariant submanifold in $\mc E$}

Let $h\in Met_{\mb Z_p}(S^3)$ be a $\langle\mk r_{p,q}\rangle$-invariant smooth metric on $S^3$ with $[h]\in \mc C_{\mb Z_p}$ satisfying $Y(S^3, [h])>0$. 
As noted at the beginning of this section, we may assume that $h$ has positive scalar curvature. 
By \cite{marques2012deforming} and \cite{chow2024equivariant}, there exists a path 
\[\{h(t)\}_{t\in[0,1]}\subset Met_{\mb Z_p}(S^3)\qquad\text{with}\qquad h(0)=h_0\quad\text{and}\quad h(1)=h\]
so that every $(S^3,h(t))$ has positive scalar curvature $R_{h(t)}>0$. 
Using the linear parameterization $\mc A$ of $\mc C$ in Lemma \ref{Lem: parameterization of boundary conformal classes}, the path $A(t):=F([h(t)])$ in $\mc A$ from $0=F([h_0])$ to $A(1)=F([h])$ is lying in an open set 
\begin{align}\label{Eq: open set of positive yamabe constant}
    \mc P:=\{A\in\mc A: Y(S^3, F^{-1}(A)) > \epsilon_0 \}
\end{align}
for some $\epsilon_0>0$. 

In this subsection, we construct a finite dimensional $\mb Z_p$-invariant topological manifold $\mc Z$ in $\mc E$ so that $\Pi(\mc Z)$ forms a connected finite dimensional smooth manifold in $\mc C$ that contains $[h]$ and the round conformal class $[h_0]$. 


\begin{proposition}\label{Prop: degree on a finite submanifold}
    There exists a finite dimensional $\langle\mk r_{p,q}\rangle$-invariant linear subspace $W\subset \mc A$ and a connected $\langle\mk r_{p,q}\rangle$-invariant open set $\Omega$ in $W\cap \mc P$ so that 
    \begin{itemize}
        \item[(i)] $\Omega_{\mb Z_p}:=\Omega\cap\mc A_{\mb Z_p}$ is a path connected submanifold containing $A(0)=0$ and $A(1)$;
        \item[(ii)] $\mc Z:=\wha \Pi^{-1}(\Omega)\subset \mc E$ is a $\langle \mk r_{p,q}\rangle$-invariant topological manifold of finite dimension $\dim(\mc Z)=\dim(W)$, and the $\langle \mk r_{p,q}\rangle$-fixed-point set $\mc Z_{\mb Z_p}:=\mc Z\cap\mc E_{\mb Z_p}\subset \wha\Pi^{-1}(\Omega_{\mb Z_p})$ is a topological submanifold of dimension $\dim(\Omega_{\mb Z_p})$; 
        \item[(iii)] the restricted maps $\wha \Pi|_\mc Z: \mc Z\to \Omega$ and $\wha\Pi|_{\mc Z_{\mb Z_p}}: \mc Z_{\mb Z_p}\to\Omega_{\mb Z_p}$ are proper maps; 
        \item[(iv)] $(\wha \Pi|_\mc Z)^{-1}(0)= (\Pi|_\mc Z)^{-1}([h_0])$ contains only one element $z_0:=[g_+^0]_{PE}$ so that $\wha \Pi|_\mc Z$ is a homeomorphism near $z_0$, where $g_+^0$ is isometric to the hyperbolic metric on $B^4$. 
    \end{itemize} 
\end{proposition}

\begin{proof}
    Consider the path $\{A(t)\}_{t\in [0,1]}\subset \mc P\cap\mc A_{\mb Z_p}$. 
    Let $N>1$ be a sufficiently large integer to be determined later. 
    For $1\leq j\leq N$, take the vertices $A_0:=A(0)=0$, $A_j:=A(\frac{j}{N})$, and the segment $[A_{j-1},A_j]$ in the linear space $\mc A$. 
    Since $\{A_j\}_{j=0}^N\subset \mc A_{\mb Z_p}$, the polyline $K:=\cup_{j=1}^N[A_{j-1},A_j]$ is also contained in the linear subspace $\mc A_{\mb Z_p}$. 
    Define 
    \[W_0:=\Span \{A_1,A_2,\dots,A_N\}\subset \mc A_{\mb Z_p}\subset \mc A\]
    to be a finite dimensional linear subspace containing $K$. 

    By taking $N>1$ sufficiently large, the compact polyline $K$ approximates the continuous path $A(t)$ and is contained in the open set $\mc P$ defined in \eqref{Eq: open set of positive yamabe constant}. 
    Let 
    \[\mc K:=\wha\Pi^{-1}(K)\subset \mc E\qquad {\rm and}\qquad \mc K_{\mb Z_p} := \mc K\cap \mc E_{\mb Z_p}. \] 
    Note that $K\subset \mc A$ is compact in the $C^\infty$ topology. 
    It then follows from the compactness theorem (Theorem \ref{Thm: compactness}) and the boundary regularity result \cite{chrsciel2005BoundaryRegularity} that $\mc K$ is a compact set in $\mc E$ (under the $C^\infty$ topology). 
    Since $\mc K_{\mb Z_p}$ is closed in $\mc K$, we also have the compactness of $\mc K_{\mb Z_p}$.

    \medskip
    \noindent{\bf Step 1.} {\it Local transversality in the coordinates.}

    Given any $z\in \mc K$, take a representative $g_+\in z$ and the $\mb Z_p$-invariant boundary metric $\gamma\in \Pi(z)$ with unit density (cf. \eqref{Eq: unit-density representative}), i.e. $d\vol_\gamma=d\vol_{h_0}$. 
    For $m\geq 6$ and $\alpha,\delta\in (0,1)$, let $\mc O_{A_z}\subset \mc A$ and $U_z\subset S^{m,\alpha}_\delta$ be small open neighborhoods of $A_z$ and $0$ respectively, where 
    \[A_z:=F(\Pi(z))\in \mc A.\]  
    Denote by 
    \[G_z: \mc O_{A_z}\times U_z\to S^{m-2,\alpha}_\delta\] 
    the Einstein operator with fixed gauge defined in \eqref{Eq: Einstein operator with gauge} with respect to $g_+\in z$ and $A_z$. 
    Then by Proposition \ref{Prop: gauge and chart}, there exists an open neighborhood $\mc U_z\subset \mc E$ of $z$ so that $\mc U_z$ can be parameterized by the solutions of $G_z=0$ in $\mc O_{A_z}\times U_z$ via the map
    \[ \Psi_z:  G_z^{-1}(0)\subset \mc O_{A_z}\times U_z \to \mc U_z,\qquad \Psi_z(A,u):=[T_{z,A}+u]_{PE},\]
    where $T_{z,A}$ is defined as in \eqref{Eq: boundary deformation} with $A_z$ in place of $A_0$. 
    For any $z'\in \mc U_z$, define
    \[(A_{z'},u_{z'}) :=\Psi_z^{-1}(z')\in \mc O_{A_z}\times U_z.\]  
    Then, for any linear subspace $V\subset \mc A$, define the operator
    \begin{align*}
        H_V(z') : V \oplus S^{m,\alpha}_\delta(B^4) &\longrightarrow S^{m-2,\alpha}_\delta(B^4),
        \\
        (v,\xi) \qquad&\longmapsto (D_AG_z)_{(A_{z'},u_{z'})}(v) + (D_uG_z)_{(A_{z'},u_{z'})} (\xi) , 
    \end{align*}
    which depends on $z'$ continuously.

    \begin{claim}\label{Claim: local transversal}
        Using the above notations, there is a finite dimensional linear subspace $V_z\subset \mc A$ and an open neighborhood $\mc U_z'\subset \mc U_z$ of $z$ so that for any $z'\in \mc U_z'$, $H_{V_z}(z')$ is surjective 
        \begin{align}\label{Eq: surjectivity of tangent map}
            \im H_{V_z}(z')  = (D_AG_{z})_{(A_{z'},u_{z'})}(V_z) + \im (D_uG_z)_{(A_{z'},u_{z'})} =S^{m-2,\alpha}_\delta(B^4),  
        \end{align}
        and $(D_uG_{z})_{(A_{z'},u_{z'})}$ is a Fredholm map of Fredholm index $0$.
    \end{claim}
    \begin{proof}
        By Proposition \ref{Prop: linearization theory}(i), the operator $L_{g_+}=\frac{1}{2}\Delta_L+3: S^{m,\alpha}_\delta(B^4)\to S^{m-2,\alpha}_\delta(B^4)$ under the metric $g_+$ is a Fredholm operator of Fredholm index $0$. 
        Hence, we can take a finite basis $\{\kappa^i\}_{i=1}^{k_z}$ of $\mc N_z:=\ker L_{g_+}$. 
        For each $1\leq i\leq k_z$, let $\kappa_{(3)}^i$ be the smooth $\gamma$-trace-free symmetric $2$-tensor on $S^3$ given by the expansion of $\kappa^i$ in Proposition \ref{Prop: linearization theory}(ii). 
    
        In particular, each $\kappa_{(3)}^i$ can be regarded as a tangent vector in $T_{[\gamma]}\mc C$. 
        Indeed, using the unique unit-density representatives in \eqref{Eq: unit-density representative}, we can identify $\mc C$ with the space of unit-density smooth metrics on $S^3$, and use $\hat{\gamma}_t^i:= (d\vol_{h_0}/d\vol_{\gamma+t\kappa_{(3)}^i})^{2/3} \cdot (\gamma+t\kappa_{(3)}^i)$ to represent the curve $[\gamma+t\kappa_{(3)}^i]$ near $[\gamma]$ in $\mc C$. 
        Combining $\frac{d}{dt}\big|_{t=0} d\vol_{\gamma + t\kappa_{(3)}^i}  = \frac{1}{2}\tr_{\gamma}(\kappa_{(3)}^i) d\vol_{\gamma}$ with $d\vol_{\gamma}=d\vol_{h_0}$ and $\tr_\gamma\kappa_{(3)}^i=0$, we have $\frac{d}{dt}\big|_{t=0} \hat{\gamma}_t^i= \kappa_{(3)}^i-\frac{1}{3}(\tr_{\gamma}\kappa_{(3)}^i)\gamma = \kappa_{(3)}^i$, which represents the tangent vector $\frac{d}{dt}\big|_{t=0}[\gamma+t\kappa_{(3)}^i]\in T_{[\gamma]}\mc C$. 
    
        Using the linear parameterization $F:\mc C\to\mc A$ (Lemma \ref{Lem: parameterization of boundary conformal classes}), we can also define
        \[A^i_z:=DF_{[\gamma]}(\kappa_{(3)}^i) \in \mc A,\qquad 1\leq i\leq k_z,\]
        which generates a finite dimensional linear subspace
        \[ V_z:=\Span\{A^1_z,\dots, A^{k_z}_z\}\subset \mc A. \]
        Then for $s=(s_1,\dots,s_{k_z})\in \mb R^{k_z}$, we have a smooth family in $\mc C$ defined by
        \[ [\gamma_s]:= F^{-1}\left(A_z + \sum_{i=1}^{k_z} s_i A^i_z \right) \in \mc C,\]
        where $\gamma_s$ is chosen to be the unit density representative \eqref{Eq: unit-density representative}. 
        In particular, by the construction of $F$ in Lemma \ref{Lem: parameterization of boundary conformal classes}, 
        \[ J_{A_z}(A^i_z)= \left.\frac{\bd \gamma_s}{\bd s_i}\right|_{s=0} = DF^{-1}_{A_z}(A^i_z)= \kappa_{(3)}^i , \]
        where the map $J_{A_z}$ is defined as in \eqref{Eq: boundary realization}. 
    
        Consider the  operator $G: \mb B^{k_z}_r(0)\times U_z\subset \mb R^{k_z}\times S^{m,\alpha}_\delta \to S^{m-2,\alpha}_\delta $ defined by 
        \[ G(s,u):= G_z\left(A_z+\sum_{i=1}^{k_z}s_i A^i_z, u \right) .\]
        By Proposition \ref{Prop: linearization theory}(iii), for any $v=(v_1,\dots,v_{k_z})\in\mb R^{k_z}$,
        \[ \int_{B^4} \langle D_sG_{(0,0)}(v),\kappa^i\rangle_{g_+} \, dV_{g_+}= -\frac{3}{2}\int_{S^3}\langle  \sum_{j=1}^{k_z} v_j\kappa_{(3)}^j , \kappa_{(3)}^i \rangle_{\gamma} \,dA_\gamma .  \]
        Meanwhile, by (ii) and (iv) in Proposition \ref{Prop: linearization theory}, the map $\kappa\mapsto \kappa_{(3)}$ is a linear injective map. 
        Therefore, the projection of $D_sG_{(0,0)}(T_0\mb R^{k_z})$ to $\mc N_z$ (with respect to the $L^2(B^4,g_+)$ inner product) is surjective, which together with $\im D_uG_{(0,0)}=\mc N_z^\perp$ in Proposition \ref{Prop: linearization theory}(i) implies that $DG_{(0,0)}$ is surjective, i.e. 
        \begin{align*}
            S^{m-2,\alpha}_\delta(B^4) = \im(DG_{(0,0)}) &= \im(D_sG_{(0,0)}) + \im(D_uG_{(0,0)}) 
            \\
            &= (D_AG_z)_{(A_z,0)}(V_z) + \im (D_uG_z)_{(A_z,0)} = \im H_{V_z}(z).
        \end{align*}
        Hence, the operator $H_{V_z}(z)$ is surjective. 
        Additionally, $(D_uG_z)_{(A_z,0)}=L_{g_+}$ is Fredholm of index $0$ (Proposition \ref{Prop: linearization theory}(i)). 
        Finally, since the surjectivity of $H_{V_z}(z')$ and the Fredholm property of $(D_uG_z)_{(A_{z'},u_{z'})}$ are open conditions, we can take a sufficiently small open neighborhood $\mc U_z'\subset \mc U_z$ of $z$ satisfying the desired result. 
    \end{proof}

    \medskip
    \noindent{\bf Step 2.} {\it Transversal extension of $W_0$ and the topological manifold $\mc M^W$ in $\mc E$.}

    For any $z\in\mc K$, let $A_z:=F(\Pi(z))\in K$,  $G_z:\mc O_{A_z}\times U_z\to S^{m-2,\alpha}_\delta$, $\Psi_z:G_z^{-1}(0)\to \mc U_{z}$, $V_z\subset \mc A$, and $\mc U_z'\subset \mc U_z$ be given as in {\bf Step 1} and Claim \ref{Claim: local transversal}. 
    If $z\in \mc K_{\mb Z_p}$, then we further shrink $\mc U_z'$ to be $\mb Z_p$-invariant in $\mc E$. 
    If $z\notin \mc K_{\mb Z_p}$, then we shrink $\mc U_z'$ so that $\mc U_z'\cap\mc E_{\mb Z_p}=\emptyset$. 
    
    By the compactness of $\mc K$ and $\mc K_{\mb Z_p}$, we obtain finitely many elements $z_1,\dots,z_k\in\mc K$ with the associated open sets $\{\mc U_{z_i}'\}_{i=1}^k$ and the finite dimensional linear subspaces $\{V_{z_i}\}_{i=1}^k$ of $\mc A$ so that 
    \[ \{z_i\}_{i=1}^{k_0}\subset \mc K_{\mb Z_p}, \qquad\mc K_{\mb Z_p}\subset \cup_{i=1}^{k_0} \mc U_{z_i}', \qquad\text{and}\qquad \mc K \subset \cup_{i=1}^{k} \mc U_{z_i}' ,\]
    where $1\leq k_0\leq k$. 
    Define
    \[ W:= W_0+ \sum_{j=1}^{p}\sum_{i=1}^k \mk r_{p,q}^j\cdot V_{z_i} \subset \mc A  \]
    as a $\langle\mk r_{p,q} \rangle$-invariant finite dimensional linear subspace of $\mc A$. 
    Then we have 
    \begin{itemize}
        \item $K\subset W_{\mb Z_p}:= W\cap \mc A_{\mb Z_p}$;
        \item $H_{W}(z')$ is surjective by \eqref{Eq: surjectivity of tangent map} for any $z'\in \mc U_{z_i}'$ and $1\leq i\leq k$. 
    \end{itemize}
    For each $1\leq i\leq k$, let $\mc B_{z_i}\subset \mc O_{A_{z_i}} \times U_{z_i}$ be an open neighborhood of $(A_{z_i}, 0)$ in $\mc A\times S^{m,\alpha}_\delta$ so that $G_{z_i}^{-1}(0)\cap \mc B_{z_i} $ is homeomorphic to $\mc U_{z_i}'$ via the map $\Psi_{z_i}$ (cf. Proposition \ref{Prop: gauge and chart}). 
    In particular, 
    \[ \mc B_{z_i}^W:= \mc B_{z_i}\cap (W\times S^{m,\alpha}_\delta) \]
    is an open neighborhood of $(A_{z_i}, 0)$ in $W\times S^{m,\alpha}_\delta$ so that 
    \begin{align}\label{Eq: local chart on Pi-1(W)}
        \Psi_{z_i}: G_{z_i}^{-1}(0)\cap \mc B_{z_i}^W \longrightarrow  \mc M_{z_i}^W:=\mc U_{z_i}' \cap \wha\Pi^{-1}(W)
    \end{align}
    is a homeomorphism. 
    Moreover, if $z_i\in\mc K_{\mb Z_p}$, then we can further choose $\mc B_{z_i}$ and $\mc B_{z_i}^W$ to be invariant under the $\mb Z_p$-action defined in Remark \ref{Rem: equivariant filling has equivariant deformation} so that $\Psi_{z_i}$ is $\mb Z_p$-equivariant. 
    
    Next, on the smooth Banach manifold $\mc B_{z_i}^W$ (open in $W\times S^{m,\alpha}_\delta$), the operator 
    \[\mc G_{z_i}:= G_{z_i}|_{\mc B_{z_i}^W}\] 
    has a surjective tangent map $(D\mc G_{z_i})_{(A_{z'},u_{z'})}=H_{W}(z')$ at any $(A_{z'},u_{z'})\in \mc G_{z_i}^{-1}(0)$. 
    Additionally, note that $(D\mc G_{z_i})_{(A_{z'},u_{z'})}=D_1+D_2$, where 
    \begin{itemize}
        \item $D_1: W\times S^{m,\alpha}_\delta\to S^{m-2,\alpha}_\delta$, $ (v,\xi)\mapsto (D_u\mc G_{z_i})_{(A_{z'},u_{z'})}(\xi)$, is a Fredholm map of index $\ind(D_1)=\dim (W)+\ind((D_u\mc G_{z_i})_{(A_{z'},u_{z'})}) = \dim (W)$, 
        \item $D_2: W\times S^{m,\alpha}_\delta\to S^{m-2,\alpha}_\delta$, $(v,\xi)\mapsto (D_A\mc G_{z_i})_{(A_{z'},u_{z'})}(v)$, is an operator of finite rank. 
    \end{itemize}
    Therefore, $(D\mc G_{z_i})_{(A_{z'},u_{z'})}$ is a Fredholm operator of index $\dim(W)$. 
    Consequently, $0$ is a regular value of $\mc G_{z_i}$ (cf. \cite{abraham1988manifolds}*{Definition 3.5.3}). 
    The implicit function theorem on Banach manifolds (cf. \cite{abraham1988manifolds}*{Theorem 3.5.4}) shows that $\mc G_{z_i}^{-1}(0)$ is a smooth submanifold in $\mc B_{z_i}^W$ of dimension $\dim(W)$. 
    
    By \eqref{Eq: local chart on Pi-1(W)}, each $\mc M_{z_i}^W$ is homeomorphic to the smooth submanifold $\mc G_{z_i}^{-1}(0)$ in $\mc B_{z_i}^W$ of dimension $\dim (W)<\infty$ via $\Psi_{z_i}$. 
    Hence, 
    \[\mc M^W:=\cup_{i=1}^k\mc M_{z_i}^W=\wha\Pi^{-1}(W)\cap (\cup_{i=1}^k\mc U_{z_i}')\] 
    is a topological manifold of dimension $\dim (W)$ with $\mc K\subset \mc M^W$. 
    (To apply a degree argument, we do not need the smoothness of the transition map $\Psi_{z_j}^{-1}\circ \Psi_{z_i}$.)

    \medskip
    \noindent{\bf Step 3.} {\it The topological manifold $\mc M^W_{\mb Z_p}$ in $\mc E_{\mb Z_p}$}
    
    For any $1\leq i \leq k_0$, $z_i=[g_+^i]_{PE}\in \mc K_{\mb Z_p}$. 
    By Proposition \ref{prop: action on the filling space}(ii), there exists $\phi_{i}\in \Diff(\overline{B}^4)$ so that $\phi_{i}|_{\bd\overline{B}^4}=\mk r_{p,q}$, $\phi_i^p=\Id_{\overline{B}^4}$, and $\phi_i^*g_+^i=g_+^i$. 
    It then follows from Remark \ref{Rem: equivariant filling has equivariant deformation} that the smooth local chart $\mc G_{z_i}^{-1}(0)$ in $\mc B^W_{z_i}\subset  W\times S^{m,\alpha}_\delta$ is invariant under the smooth $\mb Z_p$-action generated by the map $(A,u)\mapsto (\mk r_{p,q}\cdot A, \phi_i^*u)$, and the homeomorphism $\Psi_{z_i}: \mc G_{z_i}^{-1}(0)\to \mc M_{z_i}^W\subset \mc E$ is $\mb Z_p$-equivariant. 
    For any compact Lie group acting smoothly on a manifold, it is well known that the fixed-point set forms an embedded smooth submanifold. 
    In particular, the $\mb Z_p$-fixed-point set of $\mc G_{z_i}^{-1}(0)$ is a smooth submanifold that is homeomorphic to the $\mb Z_p$-fixed-point set $\mc M^W_{z_i}\cap \mc E_{\mb Z_p}$ via $\Psi_{z_i}$. 
    Hence, 
    \[\mc M_{\mb Z_p}^W:= \mc M^W\cap \mc E_{\mb Z_p} =\cup_{i=1}^{k_0} (\mc M_{z_i}^W\cap \mc E_{\mb Z_p})=\wha\Pi^{-1}(W)\cap (\cup_{i=1}^{k_0}\mc U_{z_i}') \cap \mc E_{\mb Z_p}\] 
    is a topological manifold with $\mc K_{\mb Z_p}\subset \mc M_{\mb Z_p}^W$. 
    
    Furthermore, since each $\phi_i$ is a $g_+^i$-isometry, we can restrict the linearization $(D\mc G_{z_i})_{(A_{z_i},0)}$ to the $\mb Z_p$-fixed-point subspaces $W_{\mb Z_p}:=W\cap \mc A_{\mb Z_p}$ and $S^{m,\alpha}_{\delta,\mb Z_p}:=\{\xi\in S^{m,\alpha}_\delta: \phi_i^*\xi=\xi\} $. 
    Then, the restricted map $(D_uG_{z_i})_{(A_{z_i},0)}|_{S^{m,\alpha}_{\delta,\mb Z_p}}: S^{m,\alpha}_{\delta,\mb Z_p}\to S^{m-2,\alpha}_{\delta,\mb Z_p}$ has kernel $\mc N^i_{\mb Z_p}=\ker(L_{g_+^i}) \cap S^{m,\alpha}_{\delta,\mb Z_p}$ and range $\im(L_{g_+^i})\cap S^{m-2,\alpha}_{\delta,\mb Z_p} = (\ker L_{g^i_+})^\perp\cap S^{m-2,\alpha}_{\delta,\mb Z_p} = (\mc N^{i}_{\mb Z_p})^{\perp} $ (cf. Proposition \ref{Prop: linearization theory}(i)), which implies that $(D_uG_{z_i})_{(A_{z_i},0)}|_{ S^{m,\alpha}_{\delta,\mb Z_p}}$ is a Fredholm map of index $0$. 
    By the same argument as in {\bf Step 2}, $(D\mc G_{z_i})_{(A_{z_i},0)}|_{W_{\mb Z_p}\times S^{m,\alpha}_{\delta,\mb Z_p}}$ is a Fredholm map of index $\dim(W_{\mb Z_p})$. 
    Therefore, $\dim (\mc M_{\mb Z_p}^W)= \dim(W_{\mb Z_p})$. 

    \medskip
    {\bf Step 4.} {\it The constructions of $\mc Z$ and $\mc Z_{\mb Z_p}$.}
    
    Next, we aim to construct a $\langle\mk r_{p,q}\rangle$-invariant open neighborhood $\Omega\subset W$ of $K$ so that $\mc Z:=\wha\Pi^{-1}(\Omega)\subset \mc M^W$, which implies that $\mc Z$ is a $\mb Z_p$-invariant open set in the topological manifold $\mc M^W$. 
    Firstly, take any inner product $\langle\cdot, \cdot\rangle$ on the finite dimensional subspace $W$, and consider the averaged $\langle \mk r_{p,q}\rangle$-invariant inner product $\langle u,v\rangle_{\mb Z_p}:=\frac{1}{p}\sum_{j=1}^p\langle \mk r_{p,q}^j\cdot u, \mk r_{p,q}^j\cdot v\rangle$ on $W$. 
    Hence, $\langle\cdot,\cdot\rangle_{\mb Z_p}$ induces a $\langle\mk r_{p,q}\rangle$-invariant norm $\|\cdot \|$ on $W$. 
    Since $K\subset W\cap \mc A_{\mb Z_p}$ is contained in the $\mb Z_p$-fixed-point set, the $r_0$-neighborhood $B_{r_0}(K)$ of $K$ in $(W,\|\cdot\|)$ is $\langle \mk r_{p,q}\rangle$-invariant.  
    Then, for $r_0>0$ small enough, 
    \[\Omega:=B_{r_0}(K)=\{A\in W: \inf_{A'\in K}\|A-A'\|<r_0 \} \]
    satisfies that 
    \[\Omega\subset \mc P\qquad\text{and}\qquad  \mc Z:=\wha\Pi^{-1}(\Omega)\subset \mc M^W.\]
    Indeed, $\Omega$ can be contained in $\mc P$ because $\mc P$ is open and $K\subset \mc P$. 
    If $\wha \Pi^{-1}(B_{r}(K))\not\subset\mc M^W$ for any $r>0$, then there exists a sequence $r_i\to 0^+$ so that $\emptyset\neq \wha \Pi^{-1}(B_{r_i}(K))\setminus \mc M^W$. 
    Take any $w_i\in \wha \Pi^{-1}(B_{r_i}(K))\setminus \mc M^W$ with $a_i:=\wha \Pi(w_i)$. 
    Since $\Omega\subset W$ is of finite dimension, we have the compactness of $\{a_i\}$ under the $C^{l+1,\alpha}$ topology for any $l\geq 6$.  
    Combining the compactness result (Theorem \ref{Thm: compactness}) with the boundary regularity \cite{chrsciel2005BoundaryRegularity}, $w_i$ converges up to a subsequence to an element $w\notin\mc M^W$ in the $C^{l,\alpha}$ topology for all $l\geq 6$ so that $\wha \Pi(w)=a:=\lim_{i\to\infty}a_i\in K$, which contradicts $\wha\Pi^{-1}(K)=\mc K\subset \mc M^W$. 
    
    Meanwhile, since the linear subspace $W_{\mb Z_p}:=W\cap\mc A_{\mb Z_p}$ admits the same norm $\|\cdot \|$, we conclude that $\Omega_{\mb Z_p}:=\Omega\cap \mc A_{\mb Z_p}=\{A\in W_{\mb Z_p}: \inf_{A'\in K}\|A-A'\|<r_0\}$ is path connected for $r_0>0$ small enough. 
    Moreover, by the $\mb Z_p$-equivariance of $\wha\Pi$ and the above constructions, $\mc Z_{\mb Z_p}:=\wha\Pi^{-1}(\Omega)\cap \mc E_{\mb Z_p}$ is an open set of $ \mc M^W_{\mb Z_p} $ that is contained in $ \wha\Pi^{-1}(\Omega_{\mb Z_p})$. 

    
    This proves (i) and (ii). 
    Moreover, any compact set of the finite dimensional $\Omega\subset W$ is compact in the $C^{l+1,\alpha}$ topology for all $l\geq 6$. 
    Hence, the properness of $\wha\Pi|_\mc Z$ and $\wha \Pi|_{\mc Z_{\mb Z_p}}$ also follows from the compactness result (Theorem \ref{Thm: compactness}) and the boundary regularity \cite{chrsciel2005BoundaryRegularity}.

    \medskip
    \noindent{\bf Step 5.} {\it Uniqueness of the hyperbolic filling.} 

    The statement in (iv) has been given by Anderson in \cite{anderson2008Einstein}. 
    For the sake of completeness, we include some details. 
    
    Given $0=F([h_0])\in\Omega$, let $z_0=[g_+^0]_{PE}\in \Pi^{-1}([h_0])\subset \mc E$. 
    Then there exists a geodesic defining function $\rho_0$ on $\overline{B}^4$ so that $\bar g_0:=\rho_0^2\cdot g_+^0$ extends to a $C^\infty$ metric on $\overline{B}^4$ with $h_0=\bar g_0|_{\bd \overline{B}^4}$. 

    Consider $X:=\overline{B}^4\setminus\{x_0\}$ with boundary $\bd X:=\bd \overline{B}^4\setminus\{x_0 \}\cong \mb R^3$, where $x_0\in \bd \overline{B}^4$. 
    Using the stereographic projection $P:\bd X\to \mb R^3$, we have $h_0|_{\bd X}=P^*((\frac{2}{1+|y|^2})^2\cdot \delta)$, where $y=(y_1,y_2,y_3)\in \mb R^3$ and $\delta$ is the Euclidean metric on $\mb R^3$. 
    Now, for any compact set $V\subset \bd X$, a small neighborhood $U$ of $V$ in $X$ is identified with $\mb B^3_R\times [0,\epsilon)$, for some $R,\epsilon>0$. 
    Then, for any $x\in U$ with coordinates $(y,t)\in \mb B^3_R\times [0,\epsilon)$, define $\rho(x):=\frac{1+|y|^2}{2}\rho_0(x)$. 
    Note that $d\rho=\frac{1+|y|^2}{2}d\rho_0$ on $U\cap\bd X$ since $\rho_0|_{\bd X}=0$. 
    Thus, $\rho$ is a defining function on $U$. 
    Additionally, $\rho^2g_+^0|_U=(\frac{1+|y|^2}{2})^2\bar g_0|_U$ is smooth and $\rho^2g_+^0|_{\bd X\cap U}= (\frac{1+|y|^2}{2})^2\cdot h_0|_{\bd X\cap U}=P^*\delta$. 
    Therefore, it follows from Theorem \ref{Thm: Liouville theorem for round filling} that $g_+^0$ is isometric to the hyperbolic metric $g_{\mb H}$ on $B^4$. 

    Let $\Psi: B^4\to B^4$ be a diffeomorphism with $g_+^0=\Psi^*g_{\mb H}$. 
    It follows from \cite{chrusciel2003mass}*{Theorem 6.1} that $\Psi$ extends smoothly to a diffeomorphism $\wha\Psi$ on the closed ball $\overline{B}^4$ and induces a conformal diffeomorphism $f$ on the conformal infinity $(\bd \overline{B}^4, h_0)$. 
    Since $h_0$ is the round metric on $S^3$, $f$ is a M{\"o}bius transformation that uniquely extends to a hyperbolic isometry $\hat f$ in $B^4$. 
    Therefore, $\Phi:=\hat f^{-1}\circ\wha\Psi$ satisfies $\Phi\in\Diff_\bd(\overline{B}^4)$ and $g_+^0=\Phi^*g_{\mb H}$. 
    Thus, $z_0=[g_{\mb H}]_{PE}\in \Pi^{-1}([h_0])$ is unique. 

    Finally, by \cite{lee2006fredholm}*{Theorem 5.9 and Proposition D}, the operator $(D_uG_{z_0})_{(0,0)}=L_{g_+^0}$ is an isomorphism. 
    Then, using the local chart $\Psi: \mc G_{z_0}^{-1}(0)\subset W\times S^{m,\alpha}_\delta \to \mc M^W_{z_0}$ at $z_0$ constructed as in {\bf Step 2}, 
    the implicit function theorem (\cite{abraham1988manifolds}*{Theorem 2.5.7}) gives a smooth map $u: O\to S^{m,\alpha}_\delta$ defined in a small neighborhood $O\subset W$ of $0=A_{z_0}=F([h_0])$ so that $\{(A,u(A)):A\in O\}$ is a neighborhood of $(0,0)$ in $\mc G_{z_0}^{-1}(0)$. 
    Thus, $A\in O\mapsto \Psi_{z_0}(A,u(A))$ is a homeomorphism onto a neighborhood of $z_0$ in $\mc Z$, whose inverse is $\wha \Pi|_{\mc Z}$. 
    This proves (iv). 
\end{proof}

\begin{theorem}\label{Thm: existence of P-E filling}
    For any $h\in Met_{\mb Z_p}(S^3)$ with $Y(S^3,[h])>0$, there exists a conformally compact Einstein metric $g_+$ on $B^4$ with conformal infinity $[h]$ and $[g_+]_{PE}\in \mc E_{\mb Z_p}$. 
\end{theorem}
\begin{proof}
    Given the finite dimensional topological manifold $\mc Z_{\mb Z_p}$ as in Proposition \ref{Prop: degree on a finite submanifold}, let $\mc Z_{\mb Z_p}^0$ be the connected component of $\mc Z_{\mb Z_p}$ containing the hyperbolic class $z_0=[g_+^0]_{PE}$. 
    By (ii) and (iii) in Proposition \ref{Prop: degree on a finite submanifold}, we have $k:=\dim(\mc Z_{\mb Z_p}^0)=\dim(\Omega_{\mb Z_p})<\infty$ and $\wha\Pi|_{\mc Z_{\mb Z_p}^0} :\mc Z_{\mb Z_p}^0\to\Omega_{\mb Z_p} $ is a proper map. 
    By \cite{dold1995algebraic}*{Definition 4.2}, for any compact set $K\subset \Omega_{\mb Z_p}$, we can define the mod $2$ degree $\deg_{2,K}(\wha \Pi|_{\mc Z_{\mb Z_p}^0})\in \{0,1\}$ of $\wha \Pi|_{\mc Z_{\mb Z_p}^0}$ over $K$ by 
    \[ (\wha \Pi|_{\mc Z_{\mb Z_p}^0})_*\left(o_{(\wha \Pi|_{\mc Z_{\mb Z_p}^0})^{-1}(K)}\right) = \deg_{2,K}(\wha \Pi|_{\mc Z_{\mb Z_p}^0}) \cdot o_{K}, \]
    where $o_K\in H_k(\Omega_{\mb Z_p}, \Omega_{\mb Z_p}\setminus K;\mb Z_2)$ and $o_{(\wha \Pi|_{\mc Z_{\mb Z_p}^0})^{-1}(K)}\in H_k(\mc Z_{\mb Z_p}^0, \mc Z_{\mb Z_p}^0\setminus (\wha \Pi|_{\mc Z_{\mb Z_p}})^{-1}(K);\mb Z_2)$ are the fundamental classes. 
    In particular, for a small compact neighborhood $K_0\subset \Omega_{\mb Z_p}$ of $0=\wha\Pi(z_0)$, (iv) in Proposition \ref{Prop: degree on a finite submanifold} implies that $\deg_{2,K_0}\wha \Pi|_{\mc Z_{\mb Z_p}^0}=1$. 
    
    Combining \cite{dold1995algebraic}*{Proposition and Definition 4.5} with the connectedness of $\Omega_{\mb Z_p}$ (cf. Proposition \ref{Prop: degree on a finite submanifold}(i)), the mod $2$ degree $\deg_{2,K}\wha \Pi|_{\mc Z_{\mb Z_p}^0}$ of $\wha \Pi|_{\mc Z_{\mb Z_p}^0}$ over $K$ is the same for all compact $K\subset \Omega_{\mb Z_p}$. 
    Therefore, by taking $K$ as a compact neighborhood of $F([h])$ in $\Omega_{\mb Z_p}$, the non-trivial mod $2$ degree implies $ (\Pi|_{\mc Z_{\mb Z_p}})^{-1}([h])\subset \Pi^{-1}([h])\cap \mc E_{\mb Z_p}$ is non-empty. 
\end{proof}

\part{Yamabe invariant of Lens spaces}
\section{Constructing the Test Function}

In this section, we use the existence of an equivariant Poincare-Einstein filling to calculate the Yamabe invariant of a Lens space. Fix any $p\geq 2$.  Our goal is to prove the following.

\begin{theorem}
\label{sigma-lens}
    The $\sigma$ invariant of the Lens space $L(p,q)$ is given by $\sigma(L(p,q)) = p^{-2/3}\sigma(S^3)$.
\end{theorem}

For convenience, we will write $L = L(p,q)$ from now on. Let $h$ be a smooth metric on $S^3$ which is invariant under the $\Z_p$-action generated by $\mk r_{p,q}$ so that $h$ induces a smooth metric $\bar h$ on $L$. If $Y(L,[\bar h])\leq 0$ then there is nothing to prove. So suppose that $Y(L,[\bar h]) > 0$. In this case, by Aubin's lifting lemma \cite{aubin1976scalar}, we also have $Y(S^3,[h]) > Y(L,[\bar h]) >  0$. Hence, by the results of the previous section, there exists an equivariant Poincare-Einstein filling for $[h]$.  More precisely, by Theorem \ref{Thm: existence of P-E filling} and Proposition \ref{prop: action on the filling space}, there exist the following:
\begin{itemize}
\item a conformally compact Einstein metric $g_+$ on $B^4$ satisfying $\ric(g_+) = -3g_+$, 
\item a diffeomorphism $\phi\in \Diff(\overline{B}^4)$ with $\phi|_{\bd\overline{B}^4}=\mk r_{p,q}$, $\phi^p=\Id_{\overline{B}^4}$, and $\phi^*g_+=g_+$,
    \item a $\mb Z_p$-invariant metric $\hat h\in [h]$,
    \item and a geodesic defining function $x$ with $x\circ \phi=x$ near $\bd \overline{B}^4$ so that $\hat h=(x^2 g_+)|_{\bd \overline{B}^4}$.
    \end{itemize}
    It is easy to see that $\phi$ must have a fixed point. 

    \begin{prop}
        The diffeomorphism $\phi$ has a fixed point $p_0$ in $B^4$. 
    \end{prop}

    \begin{proof}
        The Brouwer fixed point theorem implies that $\phi\colon \cl B^4 \to \cl B^4$ must have a fixed point $p_0$. Since $\phi$ agrees with $\mk r_{p,q}$ on $S^3$ and the action of $\mk r_{p,q}$ is free, it follows that $p_0$ must lie in the interior $B^4$.  
    \end{proof}

    The construction of \cite{li2017gap} now yields an invariant test function with controlled Yamabe quotient.  As in Section \ref{section:defn}, let $r = d_{g_+}(\cdot,p_0)$ and let $x$ be the invariant geodesic defining function. Set $u = r + \log x$. The function $u$ extends to a Lipschitz function $u_0$ on $S^3$ (Proposition \ref{Prop: Lipschitz extension of u=r+logx}).  Set $\psi = e^{-u/2}$ and $\psi_0 = e^{-u_0/2}$. Since $r$ and $x$ are both invariant under $\phi$, it follows that $\psi_0$ is invariant under the action of $\mk r_{p,q}$. 

    The original argument of \cite{li2017gap} uses the function $\psi_0$ to study the Yamabe quotient  via a detailed analysis of the geometry of the level sets of $r$. For the reader's convenience, we sketch another proof that illustrates the connection with the new idea in Sections \ref{distribution} and \ref{volume}. This adaptation is our own and was not suggested by the AI.  The overall strategy of the proof is exactly the same as in \cite{li2017gap}, but as in the proof of the Liouville theorem, working with regions rather than hypersurfaces leads to substantial simplifications. 
\begin{theorem}
        We have $Y(L,[\bar h]) \le p^{-2/3} \sigma(S^3)$. The equality is achieved if and only if $[\bar h]$ is the round conformal class on $L$. 
    \end{theorem}
    \begin{proof} 
     The metric $g_+$ has the form  
    \[
    g_+ = x^{-2}(dx^2 + h_x)
    \]
    where $h_x$ has the expansion 
    \[
    h_x = \hat h - x^2 \left(\ric_{\hat h} - \frac{R_{\hat h}}{4}\hat h\right) + O(x^3) = \hat h + O(x^2). 
    \]
    The geodesic estimates from Section \ref{section:defn} can still be carried through with the weaker asymptotic $h_x = \hat h + O(x^2)$, and yield that 
    \[
    0 \leq \bd_x u \leq Cx
    \]
    in a fixed size collar $\{(x,z): x\in (0,\eps),\ z\in S^3\}$.  Integrating this as $x\to 0$ gives the existence of $u_0$ satisfying 
    \[
    u(x,z) = u_0(z) + O(x^2). 
    \]
    We then set $\psi = e^{-u/2}$ and $\psi_0 = e^{-u_0/2}$. 

    Fix a small number $\eps >0$. As in Section \ref{volume}, for large $R$ there will be a unique number $0 < x_R(z)<\eps$ satisfying $r(x_R(z),z) = R$ for $z\in S^3$. Define 
    \[
    \Omega_R = \{(x,z):  x_R(z) < x < \eps,\ z\in S^3\}
    \]
    and note that $\Omega_R = B_{g_+}(p_0,R) \cap \{(x,z): x\in (0,\eps),\ z\in S^3\}$. The weaker expansion $h_x = \hat h + O(x^2)$ still suffices to carry through the volume calculation in Proposition \ref{volume-est} and thus 
    \[
    \lim_{R\to \infty} e^{-3R} \vol_{g_+}(\Omega_R) = \frac{1}{3}\int_{S^3} \psi_0^6 \, dv_{\hat h}. 
    \]
    Since $\Omega_R \subset B_{g_+}(p_0,R)$, it then follows that 
    \begin{align}\label{Eq: volume comparison}
         \int_{S^3} \psi_0^6\, dv_{\hat h} \le 3\cdot \lim_{R\to \infty} e^{-3R} \vol_{\mathbb H^4} (B(R)) = \frac{\omega_3}{8}
    \end{align}
    by the Bishop-Gromov volume comparison. 

    Next, we obtain a distributional inequality for $\psi_0$ on the boundary.  The starting point is again the Laplace comparison theorem 
    \[
    \lap_{g_+ }r \le 3\coth r. 
    \]
    Arguing exactly as in Section \ref{distribution}, one sets $F = e^{-r/2}$ and it then follows that 
    \[
    \left(\lap_{g_+} + \frac 5 4\right) F \ge -\frac{3F^5}{1-F^4}
    \]
    in the sense of distributions in the collar neighborhood of the boundary. By exactly the argument of Section \ref{distribution}, this is equivalent to 
    \[
    \psi_{xx} - \frac{\psi_x}{x} + b\psi_x + \lap_{h_x} \psi + \frac{b\psi}{2x} \ge - \frac{3\psi^5}{1-x^2 \psi^4}
    \]
    in the sense of distributions. Here $J$ is the density $dv_{h_x}/dv_{\hat h}$ and $b = \bd_x \log J$ and the precise meaning of this inequality is that 
    \[
    \int \psi \left(\eta_{xx} + \bd_x \left(\frac {\eta}{x}\right) - \bd_x(b\eta) + \frac{b}{2x} \eta + J \lap_{h_x}(J^{-1}\eta) \right) \, dx \, dv_{\hat h}(z) \ge - \int \frac{3\psi^5 \eta}{1-x^2\psi^4}\, dx\,dv_{\hat h}(z)
    \]
    for all smooth functions $\eta\geq0$ compactly supported in the collar $\{(x,z): x\in (0,\eps),\ z\in S^3\}$. 

    This inequality can then be tested against functions of the form $\phi(z)\eta(x)$ where $\phi\geq0$ is a smooth function on $S^3$ and $\eta\geq 0$ is a smooth compactly supported function on $(0,\eps)$. Letting \begin{gather*}
A_\phi(x) = \int_{S^3} \psi(x,z)\phi(z)\, dv_{\hat h} (z),\\
G_\phi(x) = \int_{S^3} \frac{b}{2x}\psi \phi + b \phi \bd_x \psi + \psi J \lap_{h_x}(J^{-1}\phi) + \frac{3\psi^5 \phi}{1-x^2 \psi^4}\, dv_{\hat h}(z),
\end{gather*}
this implies that 
\[
A_\phi '' - \frac{A_\phi'}{x} + G_\phi \ge 0
\]
in the sense of distributions on $(0,\eps)$. 

One then studies the behavior of $G_\phi(x)$ as $x\to 0$. The expansion 
\[
h_x =  \hat h - x^2 \left(\ric_{\hat h} - \frac{R_{\hat h}}{4}\hat h\right) + O(x^3) = \hat h + O(x^2)
\]
implies that 
\[
J \lap_{h_x}(J^{-1} \phi) = \lap_{\hat h} \phi + O(x^2)
\]
and 
\[
b = -\frac{R_{\hat h}}{4} x + O(x^3).
\]
Combining these estimates, it then follows that 
\[
G_\phi(x) = \int_{S^3} \psi_0 \lap_{\hat h} \phi -\frac 1 8 R_{\hat h} \psi_0 \phi + 3\psi_0^5 \phi \, dv_{\hat h}(z) + q(x)
\]
where $\vert q(x)\vert \le Cx$.  The same argument as Section \ref{distribution} now shows that necessarily one has 
\[
\int_{S^3} \psi_0 \lap_{\hat h} \phi -\frac 1 8 R_{\hat h} \psi_0 \phi + 3\psi_0^5 \phi \, dv_{\hat h}(z) \ge 0. 
\]
  Since $\phi\geq 0$ was arbitrary, this implies that 
  \[
  -\lap_{\hat h} \psi_0 + \frac 1 8 R_{\hat h}\psi_0 \le 3\psi_0^5
  \]
  in the sense of distributions on $(S^3,\hat h)$. 

  Testing this inequality against (smooth approximations to) $\psi_0$, it follows that 
  \[
  \int_{S^3} \vert \grad^{\hat h}\psi_0\vert^2 + \frac 1 8 R_{\hat h}\psi_0^2\, dv_{\hat h} \le  3 \int_{S^3}\psi_0^6 \, dv_{\hat h}.
  \]
  Hence, recalling the upper bound on the $L^6$ norm of $\psi_0$, it follows that 
  \[
  \frac{\int_{S^3} 8\vert \grad^{\hat h}\psi_0\vert^2 +  R_{\hat h}\psi_0^2\, dv_{\hat h}}{\left(\int_{S^3}\psi_0^6 \, dv_{\hat h} \right)^{1/3}} \le 24 \left(\int_{S^3} \psi_0^6\, dv_{\hat h}\right)^{2/3} \le \sigma(S^3). 
  \]
  On the other hand, abusing notation and using $\psi_0$ and $\hat h$ to denote the appropriate objects on the quotient as well, we have 
  \[
  \frac{\int_{S^3} 8 \vert \grad^{\hat h}\psi_0\vert^2 +  R_{\hat h}\psi_0^2\, dv_{\hat h}}{\left(\int_{S^3}\psi_0^6 \, dv_{\hat h} \right)^{1/3}} = p^{2/3} \left(\frac{\int_{L} 8 \vert \grad^{\hat h}\psi_0\vert^2 +  R_{\hat h}\psi_0^2\, dv_{\hat h}}{\left(\int_{L}\psi_0^6 \, dv_{\hat h} \right)^{1/3}}\right). 
  \]
  It now follows that 
  \[
  \frac{\int_{L} 8\vert \grad^{\hat h}\psi_0\vert^2 +  R_{\hat h}\psi_0^2\, dv_{\hat h}}{\left(\int_{L}\psi_0^6 \, dv_{\hat h} \right)^{1/3}} \le p^{-2/3} \sigma(S^3)
  \]
  which implies that $Y(L,[\bar h]) \le p^{-2/3}\sigma(S^3)$ as needed. 
  
  If $Y(L,[\bar h]) = p^{-2/3}\sigma(S^3)$, then $\int_{S^3}\psi_0^6\,dv_{\hat h}=\frac{\omega_3}{8}$ and the inequality in \eqref{Eq: volume comparison} given by the volume comparison achieves an equality. This implies that $(B^4,g_+)$ is isometric to hyperbolic space, as in Proposition \ref{Prop: rigidity}. 
    Thus, $[\hat{h}]$ is the round conformal class on $S^3$, and $[\bar h]$ is the round conformal class on the quotient space $L$. 
\end{proof}

\bibliographystyle{plain}
\bibliography{references.bib}

\end{document}